\magnification\magstep1


\font\ninerm=cmr9  \font\eightrm=cmr8  \font\sixrm=cmr6
\font\ninei=cmmi9  \font\eighti=cmmi8  \font\sixi=cmmi6
\font\ninesy=cmsy9 \font\eightsy=cmsy8 \font\sixsy=cmsy6
\font\ninebf=cmbx9 \font\eightbf=cmbx8 \font\sixbf=cmbx6
\font\nineit=cmti9 \font\eightit=cmti8 
\font\ninett=cmtt9 \font\eighttt=cmtt8 
\font\ninesl=cmsl9 \font\eightsl=cmsl8

\font\twelverm=cmr12 at 15pt
\font\twelvei=cmmi12 at 15pt
\font\twelvesy=cmsy10 at 15pt
\font\twelvebf=cmbx12 at 15pt
\font\twelveit=cmti12 at 15pt
\font\twelvett=cmtt12 at 15pt
\font\twelvesl=cmsl12 at 15pt
\font\twelvegoth=eufm10 at 15pt

\font\tengoth=eufm10  \font\ninegoth=eufm9
\font\eightgoth=eufm8 \font\sevengoth=eufm7 
\font\sixgoth=eufm6   \font\fivegoth=eufm5
\newfam\gothfam \def\goth{\fam\gothfam\tengoth} 
\textfont\gothfam=\tengoth
\scriptfont\gothfam=\sevengoth 
\scriptscriptfont\gothfam=\fivegoth

\catcode`@=11
\newskip\ttglue

\def\tenpoint{\def\rm{\fam0\tenrm}
  \textfont0=\tenrm \scriptfont0=\sevenrm
  \scriptscriptfont0\fiverm
  \textfont1=\teni \scriptfont1=\seveni
  \scriptscriptfont1\fivei 
  \textfont2=\tensy \scriptfont2=\sevensy
  \scriptscriptfont2\fivesy 
  \textfont3=\tenex \scriptfont3=\tenex
  \scriptscriptfont3\tenex 
  \textfont\itfam=\tenit\def\it{\fam\itfam\tenit}%
  \textfont\slfam=\tensl\def\sl{\fam\slfam\tensl}%
  \textfont\ttfam=\tentt\def\tt{\fam\ttfam\tentt}%
  \textfont\gothfam=\tengoth\scriptfont\gothfam=\sevengoth 
  \scriptscriptfont\gothfam=\fivegoth
  \def\goth{\fam\gothfam\tengoth}
  \textfont\bffam=\tenbf\scriptfont\bffam=\sevenbf
  \scriptscriptfont\bffam=\fivebf
  \def\bf{\fam\bffam\tenbf}%
  \tt\ttglue=.5em plus.25em minus.15em
  \normalbaselineskip=12pt \setbox\strutbox\hbox{\vrule
  height8.5pt depth3.5pt width0pt}%
  \let\big=\tenbig\normalbaselines\rm}

\def\ninepoint{\def\rm{\fam0\ninerm}
  \textfont0=\ninerm \scriptfont0=\sixrm
  \scriptscriptfont0\fiverm
  \textfont1=\ninei \scriptfont1=\sixi
  \scriptscriptfont1\fivei 
  \textfont2=\ninesy \scriptfont2=\sixsy
  \scriptscriptfont2\fivesy 
  \textfont3=\tenex \scriptfont3=\tenex
  \scriptscriptfont3\tenex 
  \textfont\itfam=\nineit\def\it{\fam\itfam\nineit}%
  \textfont\slfam=\ninesl\def\sl{\fam\slfam\ninesl}%
  \textfont\ttfam=\ninett\def\tt{\fam\ttfam\ninett}%
  \textfont\gothfam=\ninegoth\scriptfont\gothfam=\sixgoth 
  \scriptscriptfont\gothfam=\fivegoth
  \def\goth{\fam\gothfam\tengoth}
  \textfont\bffam=\ninebf\scriptfont\bffam=\sixbf
  \scriptscriptfont\bffam=\fivebf
  \def\bf{\fam\bffam\ninebf}%
  \tt\ttglue=.5em plus.25em minus.15em
  \normalbaselineskip=11pt \setbox\strutbox\hbox{\vrule
  height8pt depth3pt width0pt}%
  \let\big=\ninebig\normalbaselines\rm}

\def\eightpoint{\def\rm{\fam0\eightrm}
  \textfont0=\eightrm \scriptfont0=\sixrm
  \scriptscriptfont0\fiverm
  \textfont1=\eighti \scriptfont1=\sixi
  \scriptscriptfont1\fivei 
  \textfont2=\eightsy \scriptfont2=\sixsy
  \scriptscriptfont2\fivesy 
  \textfont3=\tenex \scriptfont3=\tenex
  \scriptscriptfont3\tenex 
  \textfont\itfam=\eightit\def\it{\fam\itfam\eightit}%
  \textfont\slfam=\eightsl\def\sl{\fam\slfam\eightsl}%
  \textfont\ttfam=\eighttt\def\tt{\fam\ttfam\eighttt}%
  \textfont\gothfam=\eightgoth\scriptfont\gothfam=\sixgoth 
  \scriptscriptfont\gothfam=\fivegoth
  \def\goth{\fam\gothfam\tengoth}
  \textfont\bffam=\eightbf\scriptfont\bffam=\sixbf
  \scriptscriptfont\bffam=\fivebf
  \def\bf{\fam\bffam\eightbf}%
  \tt\ttglue=.5em plus.25em minus.15em
  \normalbaselineskip=9pt \setbox\strutbox\hbox{\vrule
  height7pt depth2pt width0pt}%
  \let\big=\eightbig\normalbaselines\rm}

\def\twelvepoint{\def\rm{\fam0\twelverm}
  \textfont0=\twelverm\scriptfont0=\tenrm
  \scriptscriptfont0\sevenrm
  \textfont1=\twelvei\scriptfont1=\teni
  \scriptscriptfont1\seveni 
  \textfont2=\twelvesy\scriptfont2=\tensy
  \scriptscriptfont2\sevensy 
   \textfont\itfam=\twelveit\def\it{\fam\itfam\twelveit}%
  \textfont\slfam=\twelvesl\def\sl{\fam\slfam\twelvesl}%
  \textfont\ttfam=\twelvett\def\tt{\fam\ttfam\twelvett}%
  \textfont\gothfam=\twelvegoth\scriptfont\gothfam=\ninegoth 
  \scriptscriptfont\gothfam=\sevengoth
  \def\goth{\fam\gothfam\twelvegoth}
  \textfont\bffam=\twelvebf\scriptfont\bffam=\ninebf
  \scriptscriptfont\bffam=\sevenbf
  \def\bf{\fam\bffam\twelvebf}%
  \tt\ttglue=.5em plus.25em minus.15em
  \normalbaselineskip=12pt \setbox\strutbox\hbox{\vrule
  height9pt depth4pt width0pt}%
  \let\big=\twelvebig\normalbaselines\rm}

\def\tenbig#1{{\hbox{$\left#1\vbox
to8.5pt{}\right.\n@space$}}}
\def\ninebig#1{{\hbox{$\textfont0=\tenrm\textfont2=
\tensy\left#1\vbox to7.25pt{}\right.\n@space$}}}
\def\eightbig#1{{\hbox{$\textfont0=\ninerm\textfont2=
\ninesy\left#1\vbox to6.5pt{}\right.\n@space$}}}


\def\pp{{\goth p}}

\def\QQ{{\bf Q}}
\def\ZZ{{\bf Z}}

\def\RR{{\bf R}}
\def\CC{{\bf C}}

\def\abs#1{|\!|#1|\!|}
\def\zmod#1{\,\,({\rm mod}\,\,#1)}

\def\quod{\hskip 0.5em\relax}
\def\b{{\bf b}}

\def\exact#1,#2,#3{0\longrightarrow#1\longrightarrow#2
\longrightarrow#3\longrightarrow 0}

\def\today{\ifcase\month\or January \or February\or
March\or April\or May\or June\or July\or August\or
September\or October\or November\or December\fi
\space\number\day, \number\year}

\def\proclaim #1. #2\par{\medbreak
\noindent{\bf#1.\enspace}{\sl#2\par}\par\medbreak}

\magnification\magstep1\newcount\newpen\newpen=50
\def\qqlineA#1 #2 #3 {\line{\kern#1truept\vrule height0.5truept depth0.15truept 
width#2truept\hskip#3truept plus0.1truept\vrule height0.5truept depth0.05truept 
width16truept\hskip#3truept plus0.1truept\vrule height0.5truept depth0.15truept 
width#2truept\kern#1truept}}
\def\qqlineB#1 #2 #3 {\line{\kern#1truept\vrule height0.5truept depth0.15truept 
width#2truept\hskip#3truept plus0.1pt\vrule height0.5truept depth0.15truept 
width#2truept\kern#1truept}}
\def\qqlineC#1 #2 {\line{\hskip#1truept plus0.1pt\vrule height0.5truept 
depth0.15truept width#2truept\hskip#1truept plus0.1pt}}
\font\sanserifeight=cmss8 at 8truept
\font\sanserifeightfive=cmss8 scaled\magstep4
\newbox\boxtorre
\setbox\boxtorre=\vtop{\hsize=16truept\offinterlineskip\parindent=0truept
	\line{\vrule height1truept depth0.5truept width1truept\hfil\vrule 
	height1truept depth0.5truept width1truept\hfil\vrule height1truept 
	depth0.5truept width1truept\hfil\vrule height1truept depth0.5truept 
	width1truept\hfil\vrule height1truept depth0.5truept width1truept}
	\line{\vrule height1truept depth1truept width16truept}
	\vskip1.5truept
	\line{\vrule height1truept depth0truept width16truept}
	\vskip1truept
	\line{\vrule height1truept depth0truept width16truept}
	\vskip1truept
	\line{\vrule height1truept depth0truept width16truept}
	\vskip1truept
	\line{\vrule height1truept depth0truept width16truept}
	\vskip1truept
	\line{\vrule height1truept depth0truept width16truept}
	\vskip1truept
	\line{\vrule height1truept depth0truept width16truept}
	\vskip1truept
	\line{\vrule height1truept depth0truept width16truept}
	\vskip1truept
	\line{\vrule height1truept depth0truept width16truept}
	\vskip1truept
	\line{\vrule height1truept depth0truept width16truept}
	\vskip1truept
	\line{\vrule height1truept depth0truept width16truept}
	\vskip1truept}
\newbox\boxu
\setbox\boxu=\hbox{\vrule height1truept depth24truept width4truept\kern6truept
	\copy\boxtorre\kern6truept\vrule height1truept depth24truept 
	width4truept}
\newbox\boxuni
\setbox\boxuni=\vtop{\baselineskip=-1000pt\lineskip=-0.15truept
\lineskiplimit=0pt\parindent=0truept\hsize=36truept
\line{\copy\boxu\hfil}
\qqlineA 0.0 4.1 5.9 
\qqlineA 0.0 4.2 5.8 
\vskip-0.1truept
\qqlineB 0.1 4.2 27.4 
\qqlineB 0.21 4.2 27.2 
\qqlineA 0.3 4.2 5.5 
\qqlineA 0.41 4.2 5.4 
\qqlineB 0.5 4.3 26.4 
\qqlineB 0.6 4.4 26.0 
\qqlineA 0.7 4.5 4.8 
\qqlineA 0.8 4.7 4.5 
\qqlineB 1.0 4.8 24.4 
\qqlineB 1.2 4.8 24.0 
\qqlineB 1.5 4.8 23.4 
\qqlineB 1.7 4.91 22.8 
\qqlineB 1.9 5.1 22.0 
\qqlineB 2.21 5.2 21.2 
\vskip-0.1truept
\qqlineB 2.5 5.3 20.4 
\qqlineB 2.8 5.5 19.4 
\qqlineB 3.1 5.7 18.4 
\qqlineB 3.4 6.0 17.21 
\qqlineB 3.7 6.3 16.0 
\qqlineB 4.2 6.8 14.0 
\qqlineB 4.6 7.4 12.0 
\qqlineB 5.1 7.9 10.0 
\qqlineB 5.6 9.0 6.8  
\vskip-0.1truept
\qqlineC 6.1 23.8 
\qqlineC 6.7 22.6 
\qqlineC 7.4 21.21 
\qqlineC 8.0 20.0 
\qqlineC 8.6 18.8 
\vskip-0.1truept
\qqlineC 9.6 16.8
\qqlineC 10.91 14.2 
\qqlineC 12.0 12.0 
\qqlineC 13.5 9.0 
\vskip 3.1truept
\centerline{\sanserifeight TOR VERGATA}}

\newbox\boxprelim
\setbox\boxprelim\vtop to\dp\boxuni{\parindent=0truept\hsize=5truecm\eightpoint
Preliminary version\par\today\par\vfill}
\newbox\boxroma
\setbox\boxroma\vtop
to\dp\boxuni{\hsize=6truecm\parindent=0truept \vglue
.0truecm\sanserifeightfive II Universit\`a degli\par
\smallskip \hskip 1truecm Studi di Roma\par\vfill}

\def\title #1\par{
  \hbox{\copy\boxprelim\hskip 1.2in
\copy\boxuni\hskip 1cm\copy\boxroma}  
  \bigskip\bigskip\bigskip 
  \bigskip\bigskip\bigskip\bigskip\noindent
  {\twelvepoint$\!\!\!\!$\bf#1}\par
  \bigskip\bigskip\noindent
  Ren\'e Schoof
  \bigskip\eightpoint\noindent
  \vbox{
  \hbox{Dipartimento di Matematica}
  \hbox{$\hbox{2}^{\hbox{a}}$ Universit\`a di Roma 
  ``Tor Vergata"}
  \hbox{I-00133 Roma ITALY}
  \hbox{Email: \tt schoof@science.uva.nl}}
  \bigskip\bigskip
  \tenpoint
  }
\def\abstract #1\par{\eightpoint\vbox{\noindent
  {\bf Abstract.\ }#1\par}\bigskip\tenpoint
  }
\def\bibliography#1\par{\vskip0pt
  plus.3\vsize\penalty-200\vskip0pt
  plus-.3\vsize\bigskip\vskip\parskip
  \message{Bibliography}\leftline{\bf
  Bibliography}\nobreak\smallskip\noindent
  \ninepoint\frenchspacing#1}
\outer\def\beginsection#1\par{\vskip0pt
  plus.3\vsize\penalty\newpen\newpen=-50\vskip0pt
  plus-.3\vsize\bigskip\vskip\parskip
  \message{#1}\leftline{\bf#1}
  \nobreak\smallskip\noindent}

\def\EVA{{1}} 
\def\BUJNT{{2}} 
\def\BUJNTT{{3}}   
\def\BUHAB{{4}}
\def\BUWINF{{5}} 
\def\BUW{{6}}   
\def\BUDPP{{7}} 
\def\BUDU{{8}}   
\def\BUHO{{9}} 
\def\CO{{10}}   
\def\CDDO{{11}} 
\def\DOB{{12}} 
\def\GR{{13}}   
\def\MCCUR{{14}} 
\def\JU{{15}} 
\def\LLL{{16}}
\def\LEX{{17}}  
\def\LEB{{18}}   
\def\LEM{{19}} 
\def\LENSH{{20}} 
\def\LID{{21}}
\def\MAG{{22}}
\def\MAR{{23}}  
\def\PARI{{24}}
\def\SMC{{25}}  
\def\SH{{26}}  
\def\SHW{{27}} 
\def\SZ{{28}}  
\def\SZAST{{29}} 
\def\THIE{{30}} 
\def\EF{{31}}  
\def\WS{{32}}    
\def\WDS{{33}}    


\title Computing Arakelov class groups

\abstract Shanks's infrastructure algorithm and Buchmann's algorithm for computing class groups and unit
groups of rings of integers of algebraic number fields are most 
naturally viewed as computations inside Arakelov class
groups. In this paper we discuss the basic properties of Arakelov class groups and of the set of reduced
Arakelov divisors. As an application we describe Buchmann's algorithm in this context.

\beginsection 1. Introduction. 

In his 1972 Boulder paper~[\SH], Daniel Shanks observed that the quadratic forms 
in the principal cycle of reduced binary quadratic forms of positive discriminant exhibit a
group-like behavior. This was a surprising phenomenon, because the principal cycle itself constitutes
the trivial class of the class group. Shanks called this group-like structure `inside' the
neutral element of the class group the {\it infrastructure}. He
exploited it by designing an efficient algorithm to compute  regulators of real quadratic
number fields. Eight years later, Hendrik Lenstra  made
Shanks's observations more precise. He introduced a certain
topological group and provided a satisfactory framework for
Shanks's algorithm~[\LEX,~\SMC].  Both Shanks~[\SHW,~sect.1] and
[\WS,~sect.4.4] and Lenstra~[\LEX,~sect.15]  indicated that
the infrastructure ideas could  be generalized to arbitrary number fields.
This was done first by H.~Williams and his students~[\WDS] for
complex cubic fields, then by J.~Buchmann~[\BUJNT, \BUJNTT, \BUHAB] and by
 Buchmann and Williams~[\BUW].
Finally in~1988,  Buchmann~[ \BUDPP,~\BUDU] described  
an algorithm for computing the class group and
regulator of an arbitrary  number field that, under reasonable
assumptions,  has a subexponential running
time. It has been implemented in the LiDIA, MAGMA and PARI
software packages~[\LID,~\MAG,~\PARI].

In these expository notes we present a natural setting for the infrastructure phenomenon and
for Buchmann's algorithm. It is provided by Arakelov theory~[\SZ,~\SZAST,~\EF]. We show that
Buchmann's algorithm for computing the class number and regulator of a number field~$F$ has
a natural description in terms of  the {\it Arakelov class group} ${\rm Pic}_F^{0}$ of~$F$ and
the set ${\rm Red}_F$ of {\it reduced Arakelov divisors}.  
We show that Lenstra's topological group is essentially equal to 
the Arakelov class group of a real quadratic field.
We also introduce
the {\it oriented Arakelov class group}~$\widetilde{\rm Pic}_F^{0}$. 
This is a natural generalization of
${\rm Pic}_F^{0}$, useful for analyzing Buchmann's algorithm and for computing the units of the ring of
integers~$O_F$ themselves  rather than just the regulator.

The  main results of the paper are to be found in sections~7 and~12.
Here we discuss the basic properties of reduced Arakelov divisors
and Buchmann's algorithm respectively.

In section~2 we introduce the Arakelov class group of a number field~$F$. In section~3 we
study the \'etale $\RR$-algebra $F\otimes_{\QQ}\RR$. In section~4 we discuss the relations between Arakelov divisors, Hermitian line bundles and  ideal  lattices. In section~5 we define the 
oriented Arakelov class group and in section~6 we give both Arakelov class groups
a natural translation invariant Riemannian structure.
The rest of the notes is devoted to computational issues. Section~7 contains the main results.
Here we introduce {\it reduced} Arakelov divisors and prove their basic properties. In section~8,
we work out the details for quadratic number fields.  In section~9 we present explicit examples
illustrating various properties of reduced  divisors. In section~10 we discuss
the computational aspects of reduced Arakelov divisors. In section~11 we present a 
{\it deterministic} algorithm to compute the Arakelov class group. 
Finally, in
section~12 we present Buchmann's algorithm from the point of view of Arakelov theory. See~[\MAR] for
the basic properties of algebraic number fields.

I thank the Clay
Foundation for financial support during my stay at MSRI in the fall of~2000,  
Hendrik Lenstra and Sean Hallgren for several useful remarks  and YoungJu
Choie for inviting me to lecture on `infrastructure' at KIAS in June~2001.

\beginsection 2. The Arakelov class group.

In this section we introduce the Arakelov class group of a number field~$F$. This group is
analogous to the degree zero subgroup of the Picard group of a complete algebraic curve. In
order to have a good analogy with the geometric situation, we formally
`complete' the spectrum of the ring of integers~$O_F$ by
adjoining  primes at infinity. An {\it infinite}  prime
of $F$ is a field homomorphism $\sigma:F\longrightarrow\CC$,
considered up to complex conjugation. An infinite prime
$\sigma$ is called {\it real} when $\sigma(F)\subset\RR$ and {\it complex} otherwise.  We let $r_1$ and
$r_2$ denote the number of real and complex infinite primes, respectively. We have that
$r_1+2r_2=n$ where $n=[F:\QQ]$.

An {\it Arakelov divisor} is a formal finite sum
$D=\sum_{\pp}n_{\pp}{\pp}+\sum_{\sigma}x_{\sigma}\sigma$, where
${\pp}$ runs over the non-zero prime ideals of~$O_F$ and
$\sigma$ runs over the infinite primes of~$F$. The coefficients
$n_{\pp}$ are in~$\ZZ$ but the $x_{\sigma}$ can be any number in~$\RR$.
The Arakelov divisors form an additive  group, the  Arakelov divisor group ${\rm Div}_F$.
It is isomorphic to~$\mathop{\oplus}_{\pp}\ZZ\times\mathop{\oplus}_{\sigma}\RR$.
The {\it principal} Arakelov divisor associated to  an element $f\in F^*$ is the divisor
$(f)=\sum_{\pp}n_{\pp}\pp +\sum_{\sigma}x_{\sigma}\sigma$ with $n_{\pp}={\rm ord}_{\pp}(f)$ and
$x_{\sigma}(f)=-{\rm log}|\sigma(f)|$. The principal
Arakelov divisors form a subgroup of~${\rm Div}_F$.

 Since it is analogous to the Picard group of an algebraic curve,
 the quotient of ${\rm Div}_F$ by its
subgroup of principal Arakelov divisors is denoted by~${\rm Pic}_F$. A
principal  Arakelov divisor
$(f)$ is trivial if and only if
$f$ is a unit of $O_F$ all of whose conjugates have absolute value equal to~1. It
follows that $(f)$ is trivial if and only if $f$ is contained in the group of roots of
unity~$\mu_F$. Therefore there is an exact sequence
$$
0\longrightarrow\mu_F\longrightarrow F^* \longrightarrow{\rm Div}_F\longrightarrow{\rm
Pic}_F\longrightarrow 0.
$$
We call $I=\prod_{\pp}\pp^{-n_{\pp}}$ the {\it ideal associated} to an
Arakelov divisor~$D=\sum_{\pp}n_{\pp}\pp +\sum_{\sigma}x_{\sigma}\sigma$. 
The ideal associated to the zero Arakelov divisor is the ring of integers~$O_F$.
The ideal associated to a principal Arakelov divisor $(f)$ is
the principal ideal~$f^{-1}O_F$. Here and in the rest of the paper we often call fractional
ideals simply `ideals'. If we want to emphasize that an ideal is integral, we call it an $O_F$-ideal.

The map that sends  a divisor $D$ to its associated ideal $I$ is a homomorphism from ${\rm Div}_F$
to the group of fractional ideals ${\rm Id}_F$ of~$F$. 
 Its kernel is the group~$\mathop{\oplus}_{\sigma}\RR$
of divisors supported in the infinite primes.
 We have the following commutative
diagram, the  rows and columns of which are exact.
$$
\def\normalbaselines{\baselineskip20pt\lineskip3pt 
\lineskiplimit3pt}
\matrix{&&0&&0&&0&&\cr
&&\Big\downarrow&&\Big\downarrow&&\Big\downarrow&&\cr
0&\longrightarrow&O_F^*/\mu_F&\longrightarrow
&F^*/\mu_F&\longrightarrow&{\rm Pid}_F&\longrightarrow &0\cr
&&\Big\downarrow&&\Big\downarrow&&\Big\downarrow&&\cr
0&\longrightarrow&\mathop{\oplus}\limits_{\sigma}\RR&\longrightarrow
&{\rm Div}_F&\longrightarrow&{\rm Id}_F&\longrightarrow &0\cr
&&\Big\downarrow&&\Big\downarrow&&\Big\downarrow&&\cr
0&\longrightarrow&T&\longrightarrow
&{\rm Pic}_F&\longrightarrow&Cl_F&\longrightarrow &0\cr
&&\Big\downarrow&&\Big\downarrow&&\Big\downarrow&&\cr
&&0&&0&&0&&\cr}
$$
Here
${\rm Pid}_F$ denotes the group of principal ideals of~$F$. The map $F^*/\mu_F\longrightarrow
{\rm Div}_F$ induces a homomorphism from $O_F^*/\mu_F$ to~$\mathop{\oplus}_{\sigma}\RR$. 
 This homomorphism is given by $\varepsilon\mapsto(-{\rm log}|\sigma(\varepsilon)|)_{\sigma}$
 and its cokernel is denoted by~$T$.
 
The {\it norm} $N(\pp)$ of a non-zero prime ideal $\pp$ of~$O_F$ is the order of its residue
field~$O_F/\pp$. The {\it degree} ${\rm deg}(\pp)$ of~$\pp$  is defined
as~${\rm log}(N(\pp))$. The degree of an infinite prime
$\sigma$ is equal to~1 or~2 depending on whether $\sigma$ is real or complex. The degree extends by
linearity to a surjective homomorphism ${\rm deg}:{\rm Div}_F\longrightarrow\RR$. The
{\it norm} $N(D)$ of a divisor $D$ is defined as~$N(D)=e^{{\rm deg}(D)}$. The divisors of degree~0
form a subgroup ${\rm Div}^0_F$ of~${\rm Div}_F$. By the product formula, ${\rm Div}^0_F$ contains
the principal Arakelov divisors. 
 
 \medskip\noindent{\bf Definition 2.1.}  Let $F$ be a number field. The {\it Arakelov class group}
${\rm Pic}_F^0$  of~$F$ is the quotient of ${\rm Div}^0_F$ by its subgroup of principal divisors.
 
 \medskip
The   degree map ${\rm deg}:{\rm Div}_F\longrightarrow\RR$ factors through ${\rm Pic}_F$  and  the Arakelov class group is the kernel of the  induced homomorphism ${\rm deg}:{\rm Pic}_F\longrightarrow\RR$. We let $(\mathop{\oplus}_{\sigma}\RR)^0$ denote the subgroup of
divisors in $\mathop{\oplus}_{\sigma}\RR$ that have degree zero and $T^0$ the cokernel of the homomorphism
$O_F^*\longrightarrow (\mathop{\oplus}_{\sigma}\RR)^0$.  In other words, $T^0$ is the quotient of the vector space $\{(v_{\sigma})_{\sigma}\in\mathop{\oplus}_{\sigma}\RR:\sum_{\sigma}{\rm
deg}(\sigma)v_{\sigma}=0\}$ by the group of vectors $\{({\rm
log}|\sigma(\varepsilon)|)_{\sigma}:\varepsilon\in O_F^*\}$. By Dirichlet's unit theorem it is a
compact real torus.

\proclaim Proposition 2.2. There is a natural exact sequence
$$
0\quad\longrightarrow\quad T^0\quad \longrightarrow\quad
{\rm Pic}_F^0\quad\longrightarrow\quad Cl_F\quad\longrightarrow\quad 0.
$$

\noindent{\bf Proof.} Since $F$ has at
least one infinite prime, the composite map ${\rm Div}^0_F\hookrightarrow {\rm
Div}_F\longrightarrow{\rm Id}_F$ is still surjective. The result now follows by 
replacing the groups ${\rm Div}_F$, ${\rm Pic}_F$,  $T$ and $\oplus_{\sigma}\RR$ in the diagram
above by their degree~0 subgroups.  

\medskip
The group $T^0$ is the connected component of the identity of the topological group~${\rm Pic}^0_F$. It
follows that
${\rm Pic}^0_F$, being an extension of the finite class group by~$T^0$, is a compact real Lie group of
dimension~$r_1+r_2-1$.  

\medskip
\noindent{\bf Definition 2.3.}\ The  natural homomorphism ${\rm Div}^0_F\longrightarrow{\rm Id}_F$ admits a
 section 
$$d:{\rm Id}_F \longrightarrow{\rm Div}_F^0.
$$ 
It is given by $d(I)=D$
where $D=\sum_{\pp}n_{\pp}\pp +\sum_{\sigma}x_{\sigma}\sigma$ is the Arakelov divisor
for which we have $I=\prod_{\pp}\pp^{-n_{\pp}}$ and $x_{\sigma}=-{1\over n}{\rm log}(N(I))$ for
every infinite prime~$\sigma$. 

\proclaim Proposition 2.4. Let $\overline{d}:{\rm Id}_F\longrightarrow {\rm Pic}_F^0$ denote the homomorphism that maps $I$ to the class of the divisor $d(I)$. 
Then the sequence
$$
0\longrightarrow\quad \{f\in F^*:\hbox{all $|\sigma(f)|$ are
equal}\}/\mu_F\quad\longrightarrow\quad {\rm Id}_F\quad\mathop{\longrightarrow}\limits^{\overline{d}}\quad{\rm
Pic}^0_F
$$
is exact. Moreover, the image of~$\overline{d}$ is dense in~${\rm Pic}^0_F$.

This proposition is not used in the rest of the paper.  We do not prove it,
because  it follows immediately from Proposition~6.4 below. The kernel of~$\overline{d}$
is not a very convenient group to work with. This is one of the reasons for introducing the {\it
oriented} Arakelov divisors in section~5. 

Finally we remark that there is a natural surjective homomorphism
${\bf A}_F^*\longrightarrow {\rm Div}_F$ from the id\`ele group ${\bf A}_F^*$ to the Arakelov divisor
group. It follows that ${\rm Pic}_F$ is a quotient of the id\`ele class group. We do not make any use of this fact
in the rest of the paper.

\beginsection 3. Etale $\RR$-algebras.

Let $F$ be a number field of degree~$n$. In this section we study the
$\RR$-algebra~$F_{\RR}=F\otimes_{\QQ}\RR$.  

For any infinite prime~$\sigma$ of~$F$, we write
$F_{\sigma}$ for $\RR$ or $\CC$ depending on whether $\sigma$ is real or complex. The
natural map $F\longrightarrow\prod_{\sigma} F_{\sigma}$ that sends $f\in F$ to the vector
$(\sigma(f))_{\sigma}$, induces  an isomorphism $F_{\RR}=F\otimes_{\QQ} \RR\cong \prod_{\sigma}F_{\sigma}$
of
$\RR$-algebras.  Let $u\mapsto \overline{u}$ denote the canonical
conjugation of the \'etale algebra~$F_{\RR}$. In terms of the isomorphism
$F_{\RR}\cong\prod_{\sigma} F_{\sigma}$, it is simply the morphism that maps a
vector $u=(u_{\sigma})_{\sigma}$ to $\overline{u}=(\overline{u}_{\sigma})_{\sigma}$. In these terms
it is also easy to describe the set of invariants of the canonical conjugation. It is the subalgebra
$\prod_{\sigma}\RR$ of~$\prod_{\sigma}F_{\sigma}$.

For any $u\in F_{\RR}$, we define the norm $N(u)$ and trace ${\rm Tr}(u)$ of $u$ as the
determinant and trace respectively of the $n\times n$-matrix (with respect to any $\RR$-basis) 
of  the $\RR$-linear  map $F_{\RR}\longrightarrow F_{\RR}$ given by multiplication
by~$u$.  In terms of coordinates, we have for
$u=(u_{\sigma})_{\sigma}\in \prod_{\sigma}F_{\sigma}$ that ${\rm Tr}(u)=\sum_{\sigma}{\rm deg}(\sigma){\rm Re}(u_{\sigma})$ while
$N(u)=\prod_{\sigma}u_{\sigma}^{{\rm deg}(\sigma)}$. 

Being  an \'etale $\RR$-algebra, $F_{\RR}$  admits a canonical Euclidean 
structure~[\GR]. It is given by the scalar product
$$
\langle u,v\rangle={\rm Tr}(u\overline{v}),\qquad\hbox{for $u,v\in F_{\RR}$}.
$$
This scalar product has the `Hermitian' property $\langle\lambda u,v\rangle=\langle
u,\overline{\lambda} v\rangle$ for $u,v,\lambda\in F_{\RR}$.
In terms of coordinates, we have for $u=(u_{\sigma})_{\sigma}$ and
$v=(v_{\sigma})_{\sigma}$ in $F_{\RR}\cong\prod_{\sigma} F_{\sigma}$ that
$$
\langle u,v\rangle=\sum_{\sigma}{\rm deg}(\sigma){\rm Re}(u_{\sigma}\overline{v}_{\sigma}).
$$
We write $\abs{u}=\langle u,u\rangle^{1/2}$ for the {\it length} of~$u\in F_{\RR}$. For the element
$1\in F\subset F_{\RR}$ we have that $\abs{1}=\sqrt{n}$. For every
$u\in F_{\RR}$, all coordinates of the product
$u\overline{u}\in\prod_{\sigma}F_{\sigma}$ are  non-negative real numbers. We define $|u|$
to be the vector
$$
|u|=(|u_{\sigma}|)_{\sigma}
$$
in the group $\prod_{\sigma}\RR^*_+\,\subset\,F^*_{\RR}$.
Here we let~$\RR^*_+=\{x\in\RR^*:x>0\}$. We have $|u|^2=u\overline{u}$.
The map $u\mapsto|u|$ is a homomorphism. It is a section of the inclusion map~$\prod_{\sigma}\RR^*_+\,\,\subset\,\,F^*_{\RR}$.

\proclaim Proposition 3.1. Let $F$ be a number field of degree~$n$. For every
$u\in F_{\RR}$ we have that
\item{(i)} 
$$
N(u\overline{u})^{1/n}\le\hbox{${1\over n}$}{\rm Tr}(u\overline{u});
$$
\item{(ii)}
$$
|N(u)|\le n^{-n/2}\abs{u}^n.
$$
\vskip0in\noindent In either case, equality holds if and only if $u$ is contained in the subalgebra $\RR$ of~$F_{\RR}$.

\smallskip\noindent{\bf Proof.} Since all coordinates
of $u\overline{u}$ are non-negative, part~{\sl (i)} is just the
arithmetic-geometric mean inequality. The second inequality follows from {\sl (i)} and the fact that
$N(\overline{u})=N(u)$. 

This proves the proposition.

\beginsection 4. Hermitian line bundles and ideal lattices.

In this section we introduce the
Hermitian line bundles and  ideal lattices associated to Arakelov divisors 
and study some of their properties. 

Let $F$ be a number field of degree~$n$ and let $D=\sum_{\pp}n_{\pp}\pp
+\sum_{\sigma}x_{\sigma}\sigma$ be an Arakelov divisor.
By $I=\prod_{\pp}\pp^{-n_{\pp}}$ we denote the ideal associated to~$D$ in section~2 and by $u$ the unit $({\rm exp}(-x_{\sigma}))_{\sigma}\in\prod_{\sigma}\RR^*_+\subset F^*_{\RR}$. 
This leads to the following definition.

\medskip\noindent{\bf Definition 4.1.} Let $F$ be a number field. A {\it Hermitian line bundle} is a pair $(I,u)$
where $I$ is a fractional $F$-ideal and $u$ a unit of the algebra $F_{\RR}\cong\prod_{\sigma}F_{\sigma}$ all of whose coordinates are positive real numbers.

\medskip
As we explained above, to every Arakelov divisor $D$ there corresponds a Hermitian line bundle~$(I,u)$. This correspondence is bijective and we will often identify the two notions.
The  zero Arakelov divisor corresponds to the trivial bundle $(O_F,1)$.
A principal Arakelov divisor
$(f)$ corresponds to the Hermitian line bundle $(f^{-1}O_F,|f|)$ and the divisor $d(I)$  associated to 
a fractional ideal $I$ at the end of section~2,  corresponds to the pair~$(I,N(I)^{-1/n})$. 
Note that $N(I)^{-1/n}$ is contained in the `diagonal' subgroup~$\RR^*_+$
of~$\prod_{\sigma}\RR^*_+$.  It follows from the formulas for $N(u)$ given in the previous section
that the degree of
an Arakelov divisor $D=(I,u)$ is equal to
$-{\rm log}(|N(u)|N(I))$.

\medskip\noindent{\bf Definition 4.2.} Let $F$ be a number field.   An {\it ideal lattice} of $F$ is a projective rank~1 $O_F$-module $L$ equipped with a real valued positive definite scalar product on $L\otimes_{\ZZ}\RR$ satisfying
$\langle\lambda x,y\rangle=\langle x,\overline{\lambda} y\rangle$ for $x,y\in L\otimes_{\ZZ}\RR$
and $\lambda\in F_{\RR}$. Two  ideal lattices $L$,
$L'$ are called {\it isometric} if there is an $O_F$-isomorphism $L\cong L'$ that is 
compatible with the scalar products on $L\otimes_{\ZZ}\RR$ and~$L'\otimes_{\ZZ}\RR$.

\medskip\noindent Here $\lambda\mapsto\overline{\lambda}$  is the
canonical algebra involution of the \'etale  $\RR$-algebra $F_{\RR}$ introduced in section~3.
Note that it need not preserve~$F$. Note also that $L\otimes_{\ZZ}\RR$ has the structure of an $F_{\RR}$-module. See~[\EVA,~\GR] for more on  ideal lattices.
There is a natural way to associate an ideal lattice to an Arakelov divisor~$D$. It is
most naturally expressed in terms of the Hermitian line bundle $(I,u)$ associated to~$D$. The $O_F$-module $I$ is projective and of rank~1. Multiplication by~$u$ gives an $O_F$-isomorphism with~$uI=\{ux:x\in I\}\subset F_{\RR}$. The canonical scalar product
on~$F_{\RR}$ introduced in section~3 gives $uI$ the  structure of an ideal lattice. Alternatively, putting
$$
\abs{f}_D=\abs{uf},\qquad\hbox{for $f\in I$},
$$
we obtain a scalar product on~$I$ itself  that we extend by linearity to~$I\otimes_{\ZZ} \RR$.
 In additive notation, if $f\in I$ and  $u\in F_{\RR}^*$ is equal to
${\rm exp}((-x_{\sigma})_{\sigma})$, then $uf$ is equal 
to the vector~$(\sigma(f)e^{-x_{\sigma}})_{\sigma}\in F_{\RR}$ and we have that
$\abs{f}_D^2=\abs{uf}^2=\sum_{\sigma}{\rm deg}(\sigma)|\sigma(f)e^{-x_{\sigma}}|^2$ for $f\in I$. 
 
The ideal lattice corresponding to the zero Arakelov divisor, i.e. to the trivial bundle $(O_F,1)$, 
 is the  the ring of integers $O_F$ viewed as a subset of~$F\subset F_{\RR}$ equipped with its canonical Euclidean structure.  
The covolume of this lattice is equal to
$\sqrt{|\Delta_F|}$, where
$\Delta_F$ denotes the  discriminant of the number field~$F$. 
The covolume of the lattice associated to an arbitrary divisor $D=(I,u)$  is equal to
$$
{\rm covol}(D)=\sqrt{|\Delta_F|}\,N(I) |N(u)|
=\sqrt{|\Delta_F|}/N(D)=\sqrt{|\Delta_F|}e^{-{\rm
deg}(D)}.
$$
For any ideal $I$, the lattice  associated to the Arakelov divisor $d(I)=(I,N(I)^{-1/n})$
can be thought of as the lattice $I\subset F\subset F_{\RR}$ equipped 
 with the canonical scalar product of $F_{\RR}$, but
{\it scaled} with a factor $N(I)^{-1/n}$ so that its covolume is equal to~$\sqrt{|\Delta_F|}$.

\proclaim Proposition 4.3. Let $F$ be a number field of discriminant $\Delta_F$.
\item{(i)} The map that associates the ideal lattice $uI$ to an Arakelov divisor $D=(I,u)$, induces a bijection between the group ${\rm Pic}_F$ and the set of isometry classes of
ideal lattices.
\item{(ii)} The same map induces a bijection between the group ${\rm Pic}^0_F$ and the set of  isometry classes of
 ideal lattices of covolume~$\sqrt{|\Delta_F|}$.

\smallskip\noindent{\bf Proof.} Let $D=(I,u)$ be an Arakelov divisor and let $D'=D+(g)$ for some
$g\in F^*$. Then we have $D'=(g^{-1}I,u|g|)$ and  multiplication by $g$ induces and isomorphism $g^{-1}I\cong I$ of
$O_F$-modules. This map is also an isometry between the associated lattices since 
$$
\abs{g^{-1}f}_{D'}=\abs{u|g|g^{-1}f}=
\abs{uf}=\abs{f}_D,\quad\hbox{for all
$f\in I\otimes_{\ZZ}\RR$.}
$$
Here we use the fact that  $v=|g|g^{-1}$ satisfies $v\overline{v}=1$ and that therefore
$\abs{vh}={\rm Tr}(vh\overline{v}\overline{h})={\rm Tr}(h\overline{h})=\abs{h}$ for all~$h\in I\otimes_{\ZZ}\RR$. We conclude that the map that sends an Arakelov
divisor to its associated ideal lattice induces a well defined map from ${\rm Pic}_F$ to the set
of isometry classes of ideal lattices. This map is {\it injective}. Indeed, if
$D=(I,u)$ and
$D'=(I',u')$ give rise to isometric lattices, then there exists $g\in F^*$ so that $I'=gI$ and 
$\abs{gf}_{D'}=\abs{f}_D$ for all $f\in I\otimes_{\ZZ}\RR$. This means that 
$\abs{{u'}gf}=\abs{uf}$ for all $f\in I\otimes_{\ZZ}\RR=F_{\RR}$.
For any infinite prime~$\sigma$, we let $e_{\sigma}\in F_{\RR}$ be the idempotent for
which $\sigma(e_{\sigma})=1$ while
$\sigma'(e_{\sigma})=0$ for all~$\sigma'\not=\sigma$. Substituting $f=e_{\sigma}$,
we find that $|\sigma(g)u'_{\sigma}|=|u_{\sigma}|$ for every~$\sigma$.
It follows  that $|g|=u/u'$, implying that $D'=D+(g)$ as required.

To see that the map is
{\it surjective}, consider an ideal lattice $L$ with Hermitian scalar product $\langle\langle-,-\rangle\rangle$ on~$L\otimes_{\ZZ}\RR=F_{\RR}$. We may assume that $L$ is actually an $O_F$-ideal. The idempotent elements~$e_{\sigma}$  in $F_{\RR}\cong\prod_{\sigma}F_{\sigma}$
are invariant under the canonical involution. This implies that the 
$e_{\sigma}$ are
pairwise orthogonal because  $\langle\langle e_{\sigma},e_{\sigma'}\rangle\rangle=
\langle\langle e_{\sigma}^2,e_{\sigma'}\rangle\rangle=\langle\langle e_{\sigma},e_{\sigma}e_{\sigma'}\rangle\rangle=0$.
Therefore the real numbers $u_{\sigma}=\langle\langle e_{\sigma},e_{\sigma}\rangle\rangle^{1/2}$ determine the metric
on~$I\otimes_{\ZZ}{\RR}$. The Arakelov divisor $(L,u)$ with
$u=(u_{\sigma})_{\sigma}\in\prod_{\sigma}\RR^*_+$ is then mapped to the isometry class of~$L$.

This proves~{\sl (i)}. Part~{\sl (ii)} follows immediately from this.
\medskip
The following proposition is concerned with the lengths of the shortest non-zero vectors
in the lattices associated to Arakelov divisors.

\proclaim Proposition 4.4.  Let $F$ be a number field of degree~$n$ and let $D=(I,u)$ be an Arakelov divisor.
 Then

 \item{(i)}  for every non-zero $f$ in   $I$ we have that
$$
\abs{f}_D\ge \sqrt{n}e^{-{1\over n}{\rm deg}(D)}.
$$
Moreover, equality holds if and only if we have $D=(fO_F,\lambda|f|^{-1})$ for some $\lambda>0$. In other words,
 if and only if $D$  is
 equal to the principal Arakelov divisor~$-(f)$,  scaled by a factor~$\lambda>0$;
\item{(ii)} there exists a non-zero $f\in I$ such that $|u_{\sigma}\sigma(f)|<
\left({2\over{\pi}}\right)^{r_2/n}{\rm covol}(D)^{1/n}$ for every~$\sigma$ and hence
$$
\abs{f}_D\le\sqrt{n}\cdot\hbox{$\left({2\over{\pi}}\right)^{r_2/n}$}{\rm covol}(D)^{1/n}.
$$

Here $r_2$ is the number of complex primes of~$F$.
 
\smallskip\noindent{\bf Proof.} {\sl (i)} Let $f\in I$. By Prop.~3.1  we have that
$\abs{f}^2_D=\abs{uf}^2\ge n|N(uf)|^{2/n}$.
Since  $|N(f)|\ge N(I)$ we find that
$$
\abs{f}^2_D\ge n|N(u)N(I)|^{2/n}=ne^{-{2\over n}{\rm deg}(D)}.
$$
The last inequality follows from the fact that ${\rm deg}(D)=-\log|N(u)N(I)|$. This proves the first statement. By Prop.~1.3 equality holds if and only if all $|u_{\sigma}\sigma(f)|$
are equal to some $\lambda>0$ and if $I$ is the principal ideal generated by~$f$. This implies that
$D$ is of the form $(f^{-1}O_F,|f|^{-1}\lambda)$ as required.

To prove {\sl (ii)} consider the set 
$V=\{(y_{\sigma})_{\sigma}\in F_{\RR}:\hbox{$|y_{\sigma}|\le
 \left({2\over{\pi}}\right)^{r_2/n}{\rm covol}(D)^{1/n}$ for
all~$\sigma$}\}$. This is a bounded symmetric convex set of volume $2^{r_1}(2\pi)^{r_2}\left({2\over{\pi}}\right)^{r_2}{\rm covol}(D)=
2^n{\rm covol}(D)$. By Minkowski's Convex Body Theorem there exists
a non-zero element $f\in I$ for which $(u_{\sigma}\sigma(f))_{\sigma}\in uI\subset F_{\RR}$ is
in~$V$. This implies~{\sl (ii)}.

This proves the proposition.

\medskip

We mention the following special
case of  the proposition.

\proclaim Corollary 4.5. Let $D=(I,u)$ be an Arakelov divisor of degree~0. Then any
non-zero $f\in I$ has the property that $\abs{f}_D\ge\sqrt{n}$, with equality if and only if $D=-(f)$. On the
other hand, there exists a non-zero $f\in I$ with 
 $\abs{f}_D\le\sqrt{n}\left({2\over{\pi}}\right)^{r_2/n}\sqrt{|\Delta_F|}^{1/n}.$

\medskip
Part~{\sl (i)} of Prop.~4.4 says that the lattices associated to Arakelov divisors $D$ are rather
`nice'. They are not very skew in the sense that they do not contain any  non-zero vectors
that are extremely short with respect to~${\rm covol}(D)^{1/n}$. This property can be expressed by means of
the {\it Hermite  constant} $\gamma(D)$. The latter  is defined as the square of the length of the shortest
non-zero vector in the lattice associated to~$D$ divided by ${\rm covol}(D)^{2/n}$. The skewer the
lattice, the smaller is its Hermite constant. The constant
$\gamma(D)$ only depends on the class of~$D$ in~${\rm Pic}_F$.

\proclaim Corollary 4.6. Let $F$ be a number field of degree~$n$ and let $D=(I,u)$ be an Arakelov
divisor. Then
$$
{n\over{|\Delta_F|^{1/n}}}\,\,\,\le\,\,\,\gamma(D)\,\,\,\le\,\,\, n\hbox{$\left({2\over{\pi}}\right)^{2r_2/n}$}.
$$
The lower bound is attained if and only if $D$ is a principal divisor scaled by
some~$\lambda>0$ as in Prop.~4.4~{\sl (i)}.

 \medskip
 The function $h^0(D)={\rm log}\left(\sum_{f\in I}{\rm exp}(-\pi\abs{f}_D^2)\right)$ introduced in~[\EF] 
 and briefly discussed in section~10, is 
 related to the Hermite constant $\gamma(D)$. Indeed,
 for most Arakelov divisors $D=(I,u)$ the shortest non-zero vectors in the associated lattice
 are equal to products of a root of unity by one fixed shortest vector. Moreover, for most $D$ 
 the contributions of the zero vector and these  vectors constitute the bulk of the infinite sum $\sum_{f\in I}{\rm exp}(-\pi\abs{f}_D^2)$. Therefore, for most Arakelov divisors $D$ the quantity  $(h^0(D)-1)/w_F$ is close to ${\rm exp}(-\pi\gamma(D){\rm covol}(D)^{2/n})$. Here $w_F$ denotes the number of roots of unity in the field~$F$.

\beginsection 5. The oriented Arakelov class group.

In this section we introduce the oriented Arakelov divisor group $\widetilde{\rm
Pic}_F$ associated to a number field $F$.

In section~4 we have associated to an Arakelov divisor $D$  a Hermitain line bundle~$(I,u)$. 
 Here $I$ is an ideal and
$u$ is a unit in~$F^*_{\RR,+}=\prod_{\sigma}\RR^*_+$. An
{\it oriented}  Hermitian line bundle 
 is a pair $(I,u)$ where $I$ is an ideal and $u$ is an {\it
arbitrary} unit in
$F^*_{\RR}\cong\prod_{\sigma}F^*_{\sigma}$. The corresponding {\it oriented Arakelov divisors} 
are formal sums $\sum_{\pp}n_{\pp}\pp+\sum_{\sigma}x_{\sigma}\sigma$ with $n_{\pp}\in\ZZ$ and $x_{\sigma}\in F_{\sigma}^*$. They form a group
$\widetilde{\rm Div}_F$ and we have that
$$\widetilde{\rm Div}_F\cong{\rm Id}_F\times F^*_{\RR}\cong
\mathop{\oplus}\limits_{\pp}\ZZ\times\prod_{\sigma}F^*_{\sigma}.
$$ 
The {\it principal} oriented Arakelov divisor associated to $f\in F^*$ is simply
the oriented divisor corresponding to the oriented Hermitian bundle $(f^{-1}O_F,f)$, where the second coordinate $f$ is viewed as an element of~$F^*_{\RR}$.
The cokernel of the injective homomorphism $F^*\longrightarrow \widetilde{\rm Div}_F$
is denoted by~$\widetilde{\rm Pic}_F$.
 The inclusion ${\rm Div}_F\subset\widetilde{\rm Div}_F$  admits the natural section $\widetilde{\rm
Div}_F\longrightarrow{\rm Div}_F$ given by~$(I,u)\mapsto(I,|u|)$.
The degree ${\rm deg}({D})$  of an oriented Arakelov divisor ${D}=(I,u)$ is by
definition the degree of the `ordinary' Arakelov divisor~$(I,|u|)$. In this way principal
oriented Arakelov divisors have degree~0. 
 
 \medskip\noindent{\bf Definition 5.1.}
 The quotient of the group  $\widetilde{\rm Div}_F^0$ of
oriented Arakelov divisors of degree~0 by the subgroup of principal divisors is called the {\it oriented Arakelov class group}. It is denoted by~$\widetilde{\rm Pic}_F^0$.
 \medskip

The  commutative
diagram below has exact rows and columns. The bottom row relates the groups $\widetilde{\rm Pic}_F^0$
and
${\rm Pic}_F^0$ to one another. 
$$
\def\normalbaselines{\baselineskip20pt\lineskip3pt 
\lineskiplimit3pt}
\matrix{&&0&&0&&0&&\cr
&&\Big\downarrow&&\Big\downarrow&&\Big\downarrow&&\cr
0&\longrightarrow&\mu_F&\longrightarrow
&F^*&\longrightarrow&F^*/\mu_F&\longrightarrow &0\cr
&&\Big\downarrow&&\Big\downarrow&&\Big\downarrow&&\cr
0&\longrightarrow&\prod_{\sigma}K_{\sigma}&\longrightarrow
&\widetilde{\rm Div}_F^0&\longrightarrow&{\rm Div}_F^0&\longrightarrow &0\cr
&&\Big\downarrow&&\Big\downarrow&&\Big\downarrow&&\cr
0&\longrightarrow&(\prod_{\sigma}K_{\sigma})/\mu_F&\longrightarrow
&\widetilde{\rm Pic}^0_F&\longrightarrow&{\rm Pic}^0_F&\longrightarrow &0\cr
&&\Big\downarrow&&\Big\downarrow&&\Big\downarrow&&\cr
&&0&&0&&0&&\cr}
$$
Here $K_{\sigma}$ denotes the maximal compact subgroup of~$F_{\sigma}^*$. In other words
$K_{\sigma}=\{1,-1\}$ if $\sigma$ is real, while
$K_{\sigma}=\{z\in\CC^*:|z|=1\}$ if $\sigma$ is complex.  Since
${\rm Pic}^0_F$ and the groups $K_{\sigma}$ are compact, it follows from the exactness
of the bottom row of the diagram  that~$\widetilde{\rm Pic}^0_F$ is compact as well.

In order to see the topological structure of  $\widetilde{\rm Pic}_F^0$ better,
we  construct a second exact sequence. Let  $F^*_{\RR,\,{\rm conn}}$ 
denote the connected component of
$1\in F^*_{\RR}$. It is isomorphic to a  product of copies of $\RR^*_+$
for the real primes and $F_{\sigma}^*=\CC^*$ for the complex ones. 
It is precisely the kernel of the homomorphism $\widetilde{\rm Div}_F\longrightarrow{\rm
Id}_F\times\prod_{\sigma\,\,\rm real}\{\pm1\}$, given by mapping $D=(I,u)$ to $(I,{\rm sign}(u))$. Here
${\rm sign}(u)$ denotes the vector~$({\rm
sign}(u_{\sigma}))_{\sigma\,\,\rm real}$. 

\medskip\noindent{\bf Definition 5.2.}
By $\widetilde{T}$ we denote the
quotient of the group $F^*_{\RR,\,{\rm conn}}$ by its subgroup $O_{F,+}^*=\{\varepsilon\in
O_F^*:\hbox{$\sigma(\varepsilon)>0$ for all real $\sigma$}\}$. Taking degree zero subgroups, we put
$(F^*_{\RR,\,{\rm conn}})^0 =\{u\in F^*_{\RR,\,{\rm conn}}:N(u)=1\}$
and  $\widetilde{T}^0=(F^*_{\RR,\,{\rm conn}})^0/O_{F,+}^*$. 

\medskip

The map $\widetilde{T}\longrightarrow \widetilde{\rm Pic}_F$ given by $v\mapsto(O_F,v)$ is a well
defined homomorphism. So is the map
$\widetilde{\rm Pic}_F\longrightarrow Cl_{F,+}$ that sends the class of the divisor $(I,u)$
to the narrow ideal class of $gI$ where $g\in F^*$ is any element
for which ${\rm sign}(g)={\rm sign}(u)$.
Here the {\it narrow} ideal  class group 
$Cl_{F,+}$ is defined as the group of ideals modulo the principal ideals that are
generated by $f\in F^*_+=\{f\in F^*:\hbox{$\sigma(f)>0$ for all real $\sigma$}\}$. 
It is a finite group.

The following proposition 
says that the groups $\widetilde{T}$  and $\widetilde{T}^0$ are the connected components of identity of~$\widetilde{\rm Pic}_F$ and $\widetilde{\rm Pic}_F^0$ respectively.
It provides an analogue to Proposition~2.2.

\proclaim Proposition 5.3. Let $F$ be a number field of degree~$n$.
\item{(i)} The natural sequences
$$
0\quad\longrightarrow\quad\widetilde{T}
\quad\longrightarrow\quad\widetilde{\rm
Pic}_F\quad\longrightarrow\quad Cl_{F,+}\quad\longrightarrow \quad0
$$
and
$$
0\quad\longrightarrow\quad\widetilde{T}^0
\quad\longrightarrow\quad\widetilde{\rm
Pic}^0_F\quad\longrightarrow\quad Cl_{F,+}\quad\longrightarrow \quad0
$$
are exact.
\item{(ii)} The groups $\widetilde{T}$ and $\widetilde{T}^0$ are the connected components of identity of
$\widetilde{\rm Pic}_F$ and $\widetilde{\rm Pic}_F^0$ respectively. The  group~$\widetilde{T}$ has dimension~$n$ while $\widetilde{T}^0$ is a compact torus of dimension~$n-1$.

\smallskip\noindent{\bf Proof.} {\sl (i)}  Let $\widetilde{\rm Pid}_F$ denote   the image of  the map
 $F^*\longrightarrow \widetilde{\rm Div}_F\longrightarrow{\rm Id}_F\times\prod_{\sigma\,\,\rm real}\{\pm 1\} $. 
 This leads to the following commutative diagram with exact rows.
 \smallskip
$$
\def\normalbaselines{\baselineskip20pt\lineskip3pt 
\lineskiplimit3pt}
\matrix{0&\longrightarrow&O^*_{F,+}&\longrightarrow&F^*&\longrightarrow&\widetilde{\rm Pid}_F
&\longrightarrow&0\cr
&&\Big\downarrow&&\Big\downarrow&&\Big\downarrow&&\cr
0&\longrightarrow&F^*_{\RR,\,{\rm conn}}&\longrightarrow&\widetilde{\rm Div}_F&\longrightarrow&{\rm
Id}_F\times\prod_{\sigma\,\,\rm real}\{\pm 1\} &\longrightarrow&0.\cr}
$$
The vertical maps in this diagram are all injective. 
 An application of the snake lemma shows that the sequence of cokernels
is exact. This is  the  first exact sequence of~{\sl (i)}. Indeed,
the kernel of the surjective homomorphism ${\rm Id}_F\times\prod_{\sigma\,\,\rm
real}\{\pm1\}\longrightarrow Cl_{F,+}$ given by mapping a pair $(I,s)$ to the narrow ideal class of $gI$ where
$g\in F^*$ is any element for which ${\rm sign}(g)=s$, is precisely equal to~$\widetilde{\rm Pid}_F$. 
The second exact sequence is obtained  by taking degree zero parts.
\smallskip\noindent
{\sl (ii)} Since $Cl_{F,+}$ is finite and both groups $\widetilde{T}$ and $\widetilde{T}^0$ are connected, the first statement is clear. 
Since the  Lie group $F^*_{\RR,\,{\rm conn}}$ has dimension~$n$,
so do the groups $\widetilde{T}$ and $\widetilde{\rm Pic}_F^0$. It follows that the groups $\widetilde{T}^0$ and $\widetilde{\rm Pic}_F^0$ have dimension~$n-1$.

This proves the proposition.

\medskip

Note that the classes of two extended Arakelov divisors $(I,u)$ and $(J,v)$ are on the same
connected component of $\widetilde{\rm Pic}_F^0$ if and only if $J=gI$ for some $g\in F^*$ for
which
$u_{\sigma}v_{\sigma}g(\sigma)>0$ for each real~$\sigma$.

\medskip\noindent{\bf Definition 5.4.} An {\it embedded ideal lattice}
is  an ideal lattice $L$ together with an $O_F$-linear
isometric embedding $L\hookrightarrow F_{\RR}$.
\medskip
To every oriented Arakelov divisor ${D}=(I,u)$ we associate the  ideal lattice $uI$
together with the embedding $uI\subset F_{\RR}$. 
Two embedded ideal lattices are called isometric 
if there is an isometry of ideal lattices that commutes with the embeddings.
 We have the following analogue of Proposition~4.3.

\proclaim Proposition 5.5. Let $F$ be a number field of discriminant~$\Delta_F$. Then
\item{(i)} the map that associates to an oriented Arakelov divisor $D=(I,u)$ its 
associated embedded ideal lattice, induces a 
bijection between the  oriented Arakelov class group
$\widetilde{\rm Pic}_F$ and the set of  isometry classes of embedded ideal lattices;
\item{(ii)} the same map induces a  bijection between
$\widetilde{\rm Pic}_F^0$ and the set of isometry classes of embedded ideal lattices of covolume~$\sqrt{|\Delta_F|}$.

\smallskip\noindent{\bf Proof.}  If  two oriented Arakelov divisors
$D=(I,u)$ and $D'=(I',u')$ differ by a principal divisor $(f^{-1}O_F,f)$, then multiplication by~$f$
induces an isometry between the embedded lattices $uI$ and $u'I'$. Therefore the map in {\sl (i)} 
is well defined.
If the embedded lattices $uI$ and $u'I'$ are isometric, then this isometry is given 
by multiplication by some $x\in F_{\RR}^*$. Then $f=u^{-1}xu'$ is contained in  $F^*$ and
we have that $D-D'=(f^{-1}O_F,f)$. This shows that the map is injective. To see that the map is surjective, 
let $I$ be a fractional ideal and let $\iota:I\hookrightarrow F_{\RR}$ be an $O_F$-linear embedding.
Tensoring $I$ with~$\RR$, we obtain an $F_{\RR}$-linear isomorphism $F_{\RR}\cong I\otimes_{\ZZ}\RR\longrightarrow F_{\RR}$ which is necessarily multiplication by some~$u\in F_\RR^*$. Therefore $\iota(I)=uI$ and the oriented divisor $(I,u)$ maps to the embedded ideal lattice
$\iota:I\hookrightarrow F_{\RR}$.

This proves the proposition
\medskip
We will not use this in the rest of the paper, but note that there is a natural surjective homomorphism 
 from the id\`ele group ${\bf A}_F^*$ to the oriented Arakelov divisor group~$\widetilde{\rm Div}_F$.
It follows that the group $\widetilde{\rm Pic}_F$ is a quotient of the id\`ele class group.

\beginsection 6. Metrics on Arakelov class groups.

Let $F$ be a number field. In this section we provide the 
Arakelov class groups ${\rm Pic}_F$ and $\widetilde{\rm Pic}_F$ with translation invariant 
Riemannian structures.

By the diagram in section~2, the connected component $T$ of the group ${\rm Pic}_F$
is isomorphic to $\mathop{\oplus}_{\sigma}\RR$ modulo the closed
discrete subgroup $\Lambda=\{(\log|\sigma(\varepsilon)|)_{\sigma}: \varepsilon\in O_F^*\}$. 
Therefore the tangent space at~0 is isomorphic to
$\mathop{\oplus}_{\sigma}\RR$. Identifying  this vector space 
with the  subalgebra $\prod_{\sigma}\RR$ of $F_{\RR}=\prod_{\sigma}F_{\sigma}$, it inherits
the canonical scalar product from~$F_{\RR}$. Since this $\RR$-valued scalar product 
is positive definite, both groups  $T$ and ${\rm Pic}_F$  are in this way equipped with 
a translation invariant Riemannian structure.

For $u\in \prod_{\sigma}\RR^*_+\subset F_{\RR}^*$ we let $\log(u)$
denote the element $(\log(\sigma(u)))_{\sigma}\in\prod_{\sigma}\RR\subset F_{\RR}$.
We have that
$$
\abs{\log(u)}^2=\sum_{\sigma}{\rm deg}(\sigma)|\log(\sigma(u))|^2.
$$
\smallskip\noindent{\bf  Definition~6.1.}  For $u\in T$ we put
$$\eqalign{
\abs{u}_{\rm Pic}&= \mathop{\rm min}\limits_{{u'\in F^*_{\RR,+}}
\atop{u'\equiv u\zmod{\Lambda}}} \abs{\log(u')},\cr
&=\quad\mathop{\rm min}\limits_{\varepsilon\in O_F^*}\abs{\log(|\varepsilon| u)}.\cr}
$$
Every divisor class in~$T$ is represented by a divisor of the form $D=(O_F,u)$ for
some $u\in \mathop{\oplus}_{\sigma}\RR^*_+$.  Here $u$ is unique up 
to multiplication by units $\varepsilon\in O_F^*$. For such a divisor  class in~$T$ we define
$$
\abs{D}_{\rm Pic}=\abs{u}_{\rm Pic}.
$$
The function $\abs{u}^{}_{\rm Pic}$ on $T$ satisfies the triangle inequality. It gives rise to a distance
function  that induces the natural topology of ${\rm Pic}_F$. The distance is only defined for divisor classes $D$ and~$D'$ that lie on the same
connected component. By Prop.~2.2, the class of the difference $D-D'$ is then equal to $(O_F,u)$ for some unique $u\in T$ and we  define the {\it distance} $\abs{D-D'}_{\rm Pic}$ between $D$ and $D'$ as~$\abs{u}_{\rm Pic}^{}$.  The closed subgroups $T^0$ and ${\rm Pic}_F^0$ inherit their Riemannian structures from~${\rm Pic}_F$.

The Euclidean structures of the ideal lattices corresponding to 
Arakelov divisors and the metric on ${\rm Pic}_F$ are not
unrelated. The following proposition says that  
the difference between the Euclidean structures
of two Arakelov divisors $D,D'$ is bounded in terms of~$\abs{D-D'}_{\rm Pic}$.

\proclaim Proposition 6.2. Let $F$ be a number field and let $D=(I,u)$ and $D'=(I,u')$
be two Arakelov divisors. Then there exists a unit $\varepsilon\in O_F^*$ for
which the divisor $D''=(I,u'|\varepsilon|)$ satisfies 
$$
e^{-\abs{D-D''}_{\rm Pic}}\quad\le\quad{{\abs{x}_D}\over {\abs{x}_{D''}}}
\quad\le\quad e^{\abs{D-D''}_{\rm Pic}},\qquad\hbox{for every $x\in I$.}
$$
Note that the classes of $D'$ and $D''$ in ${\rm Pic}_F$ are the same, so that we have
$\abs{D-D'}_{\rm Pic}=\abs{D-D''}_{\rm Pic}$.

\smallskip\noindent{\bf Proof.}   Let $\varepsilon\in O_F^*$ be such that the expression
$\sum_{\sigma}{\rm deg}(\sigma){\rm log}|\sigma(\varepsilon)u'_{\sigma}/u_{\sigma}|^2$ is minimal.
Let $D''=(I,|\varepsilon| u')$. Putting  $v=u'|\varepsilon|/u$  we have as a consequence
$$\abs{D-D''}^2_{\rm Pic}=\sum_{\sigma}{\rm deg}(\sigma){\rm log}|v_{\sigma}|^2.
$$ 
For any $x\in I$ we have
$$\eqalign{
\abs{x}_{D''}^2\quad&=\quad\abs{u'\varepsilon x}^2\quad=\quad\abs{vux}^2\quad=\quad\sum_{\sigma}{\rm deg}(\sigma)|u_{\sigma}v_{\sigma}\sigma(x)|^2\quad\le\cr
&\quad\le \mathop{\rm max}\limits_{\sigma}|v_{\sigma}|^2
\sum_{\sigma}{\rm deg}(\sigma)|u_{\sigma}\sigma(x)|^2\quad=\quad 
\left(\mathop{\rm max}\limits_{\sigma}|v_{\sigma}|\right)^2
\abs{x}_{D}^2.}
$$
Since we have 
$$\log\left(\mathop{\rm max}\limits_{\sigma}|v_{\sigma}|\right)=
\mathop{\rm max}\limits_{\sigma}\log|v_{\sigma}|\le
\mathop{\rm max}\limits_{\sigma}\left|\log|v_{\sigma}|\right|\le 
\abs{D-D''}_{\rm Pic},
$$ 
the first inequality follows. The second one follows by symmetry.

\medskip
We now define a similar metric on  the {\it oriented} Arakelov class group.
By Prop.~5.3, the connected component of $\widetilde{\rm Pic}_F$ is $\widetilde{T}=F^*_{\RR,\,\rm conn}/O^*_{F,+}$. We recall that $F^*_{\RR,\,\rm conn}$ is the connected component of identity of the group~$F_{\RR}^*$. It is isomorphic to
a product of copies of~$\RR^*_+$, one for each real prime, and $\CC^*$ one for each complex prime. The group $O^*_{F,+}$ is the subgroup of $\varepsilon\in O_F^*$ for which $\sigma(\varepsilon)>0$ for every real infinite prime~$\sigma$.

The exponential homomorphism ${\rm
exp}:F_{\RR}\longrightarrow F^*_{\RR}$
is defined in terms of the usual exponential function by ${\rm exp}(u)=({\rm exp}(u_{\sigma}))_{\sigma}$ for
$u=(u_{\sigma})\in F_{\RR}\cong\prod_{\sigma} F_{\sigma}$.
The  image of the exponential function is precisely the group $F^*_{\RR,\,\rm conn}$.
The counterimage of $O^*_{F,+}$ is a discrete closed subgroup ${\Lambda}$ of $F_{\RR}$.
We have a natural isomorphism of Lie groups ${\rm exp}: F_{\RR}/{\Lambda}\mathop{\longrightarrow}\limits^{\cong}
F^*_{\RR,\,\rm conn}/O^*_{F,+}=\widetilde{T}$.  Therefore the tangent 
space of~$\widetilde{T}$ at~0 is isomorphic to~$F_{\RR}$.  The canonical scalar product on~$F_{\RR}$ provides both groups $\widetilde{T}$ and ${\widetilde{\rm Pic}}$ with a translation invariant Riemannian structure.  
 
 \medskip\noindent{\bf  Definition~6.3.}  For $u\in \widetilde{T}$ we put
 $$
\abs{u}_{\widetilde{\rm Pic}}=\mathop{\rm min}\limits_{{y\in F_{\RR}}
\atop{{\rm exp}(y)\equiv u\zmod{O^*_{F,+}}}}\abs{y}
$$
Explicitly, 
for $u\in F^*_{\RR}=\prod_{\sigma} F^*_{\sigma}$  we let $\log(u)$
denote the element $(\log(\sigma(u))_{\sigma}\in\prod_{\sigma}F_{\sigma}\subset F_{\RR}$.
Here we use the principal branch of the complex logarithm.
We have that
$$
\abs{u}_{\widetilde{\rm Pic}}\,\,=\,\,\mathop{\rm min}\limits_{\varepsilon\in O^*_{F,+}}\abs{\log(\varepsilon u)}\,\,=\,\,
\mathop{\rm min}\limits_{\varepsilon\in O^*_{F,+}}\sum_{\sigma}{\rm deg}(\sigma)|\log(\sigma(\varepsilon u))|^2.
$$
Every divisor class in~$\widetilde{T}$ can be represented by a divisor of the form $D=(O_F,u)$ for
some~$u\in F_{\RR,\,\rm conn}^*$.  Here $u$ is unique up 
to multiplication by units $\varepsilon\in O_{F,+}^*$. For any divisor $D$ of the form $(O_F,u)$
with $u\in \widetilde{T}$  we define
$$
\abs{D}_{\widetilde{\rm Pic}}=\abs{u}_{\widetilde{\rm Pic}}.
$$
The function $\abs{u}^{}_{\widetilde{\rm Pic}}$ on $\widetilde T$ satisfies the triangle inequality and  this gives rise to a distance
function that induces the natural topology on $\widetilde{\rm Pic}_F$. The distance is only defined for divisor classes $D$ and~$D'$ that lie on the same
connected component. By Prop.~5.3, the class of the difference $D-D'$ is then equal to $(O_F,u)$ for some unique $u\in \widetilde{T}$ and we  define the {\it distance} $\abs{D-D'}_{\widetilde{\rm Pic}}$ between $D$ and $D'$ as~$\abs{u}_{\widetilde{\rm Pic}}^{}$. 

The closed subgroups~$\widetilde{T}^0$ and $\widetilde{\rm Pic}_F^0$ inherit their Riemannian structures from~${\rm Pic}_F$.  We leave to  the task of proving an `oriented' 
version of Proposition~6.2 the reader.

The morphism $d:{\rm Id}_F\longrightarrow\widetilde{\rm Div}_F^0$ given by $d(I)=(I,N(I)^{-1/n})$
is a section of the natural map $\widetilde{\rm Div}_F^0\longrightarrow{\rm Id}_F$. 
The embedded ideal lattice associated to~$d(I)$ is the ideal lattice $I\subset F_{\RR}$  scaled by a factor~$N(I)^{-1/n}$. This lattice has  covolume~$\sqrt{|\Delta_F|}$.  

Next we prove an oriented version of Proposition~2.4. 
It says that the  classes of
the divisors of the form $d(I)$ are dense in~$\widetilde{\rm Pic}_F^0$
and it  implies Proposition~2.4. The exactness of
the first sequence of~[\LEX, section 9] is a special case.

\proclaim Proposition 6.4. Let $F$ be a number field of degree~$n$. Let $\overline{d}:{\rm Id}_F\longrightarrow 
\widetilde{\rm Pic}^0_F$ be the map that sends $I$ to the class of the oriented Arakelov divisor
$d(I)$ in $\widetilde{\rm Pic}^0_F$.
Then the  sequence
$$
0\longrightarrow\quad{\rm Id}_{\QQ}\quad\longrightarrow\quad{\rm Id}_{F}\quad
\mathop{\longrightarrow}\limits^{\overline{d}}\quad\widetilde{\rm Pic}^0_F
$$
is exact. The image of the map~$\overline{d}$ is dense in $\widetilde{\rm Pic}^0_F$.

\smallskip\noindent{\bf Proof.} Every ideal in ${\rm
Id}_{\QQ}$ is generated by some $f\in\QQ^*_{>0}$. Let $f\in\QQ^*_{>0}$. Then $\overline{d}$ maps the $F$-ideal $fO_F$  to the
class of the oriented Arakelov divisor $(fO_F,|N(f)|^{-1/n})$. 
Since we have $|N(f)|=|f|^n$, this divisor is equal to $(fO_F,f^{-1})$.
Therefore its image in $\widetilde{\rm Pic}^0_F$  is trivial. 
Conversely, suppose that a fractional ideal $I$ has the property that
the class of $(I,N(I)^{-1/n})$ is trivial  in~$\widetilde{\rm Pic}^0_F$. That means that
$I=fO_F$ for some $f\in F^*$ and that~$f=N(I)^{1/n}$. In other words,
$\sigma(f)=N(I)^{1/n}$ for all infinite primes~$\sigma$. This implies that all conjugates of $f$
are equal, so that
$f\in\QQ^*$. This shows that the sequence is exact.

To show that the image of $\overline{d}$ is dense, we let $0<\epsilon<1$ and
pick $D=(I,u)\in\widetilde{\rm Div}_F^0$. Note that we have $N(I)|N(u)|=1$. Consider the set
$$B=\{(v_{\sigma})_{\sigma}\in F_{\RR}: \hbox{$|v_{\sigma}-u_{\sigma}|<\epsilon|u_{\sigma}|$
for all~$\sigma$}\}.
$$
Then $B$ is a an open subset of~$F_{\RR}^*$ and all $v\in B$ have the same signature as~$u$.
Since $F$ is dense in $F_{\RR}$, there is an element $f\in B\cap F$.

The difference between $d(fI)$   and the  divisor $D$
 is equal to $(fO_F,N(fI)^{-1/n}u^{-1})$ which is
equivalent to the Arakelov divisor
$(O_F^{},v)$ where $v=N(f/u)^{-1/n}u^{-1}f\in F^*_{\RR}$.
Therefore the distance between
$D$ and $\overline{d}(f^{-1}I)$ is at most
$\abs{v}^{}_{\widetilde{\rm Pic}}$. Since 
$|{{\sigma(f)}\over{u_{\sigma}}}-1|<\epsilon$, it follows from the Taylor series expansion of the principal 
branch of the logarithm that $|{\rm
log}({{\sigma(f)}\over{u_{\sigma}}})|<{{\epsilon}\over{1-\epsilon}}$ for all~$\sigma$ and 
hence ${1\over n}|{\rm log}({{N(f)}\over{N(u)}})|<{{\epsilon}\over{1-\epsilon}}$.
It follows that we have
$$\eqalign{
\abs{v}^{}_{\widetilde{\rm Pic}}\,\,&\le\,\,\sqrt{n}\mathop{\rm max}\limits_{\sigma}|{\rm
log}(N(f/u)^{-1/n}u_{\sigma}^{-1}\sigma(f))|,\cr
&\le\,\,\sqrt{n}\left(
\hbox{${1\over n}$}|{\rm log}({{N(f)}\over{N(u)}})|+\mathop{\rm max}\limits_{\sigma}|{\rm
log}({{\sigma(f)}\over{u_{\sigma}}})|\right),\cr
&\le\,\,{2{\epsilon\sqrt{n}}\over{1-\epsilon}}.\cr}
$$ 
This implies that the image of~$\overline{d}$ is dense, 
as required.

Finally we compute the volumes of the compact Riemannian manifolds ${\rm Pic}_F^0$ and~$\widetilde{\rm Pic}_F^0$.

\proclaim Proposition 6.5.  Let $F$ be a number field of degree~$n$ and discriminant~$\Delta_F$.
Then
\item{(i)} 
$$
{\rm vol}({\rm Pic}^0_F)=
{{w_F\,\sqrt{n}}\over{2^{r_1}(2\pi\sqrt{2})^{r_2}}}\cdot|\Delta_F|^{1/2}\cdot
\mathop{\rm
Res}\limits_{s=1}\zeta_F^{}(s).
$$
\item{(ii)}
$$
{\rm vol}(\widetilde{\rm Pic}^0_F)=\sqrt{n}\cdot{|\Delta^{}_F|}^{1/2}\cdot\mathop{\rm
Res}\limits_{s=1}\zeta^{}_F(s).
$$
Here $r_1$ is the number of real primes and $r_2$ is the number of complex primes of~$F$.
By $w_F$ we denote the number of roots of unity and by $\zeta_F(s)$ the Dedekind zeta function of~$F$.

\smallskip\noindent{\bf Proof.} {\sl (i)}

The subspace $(\oplus_{\sigma}\RR)^0$ of divisors of degree~0 is the
orthogonal complement of~$1$ in the subalgebra $\prod_{\sigma}\RR$ of~$F_{\RR}$.  Using the fact that
$\abs{1}=\sqrt{n}$, one checks that the volume of~${\rm Pic}_F^0$ is equal to 
$\sqrt{n}\,2^{-r_2/2}R_F$ where $R_F$ is the regulator of~$F$.
It follows from the  exact sequence of Prop.~2.2 that the compact group
${\rm Pic}^0_F$ has volume~$\sqrt{n}\,2^{-r_2/2}h_FR_F$ where $h_F=\#Cl_F$ is the class number
of~$F$.  The formula~[\MAR]  for the residue of the zeta function in $s=1$ now easily implies~{\sl (i)}.

\noindent {\sl (ii)} Since the natural volume of
the group $K_{\sigma}$ is~2 or $2\pi\sqrt{2}$ depending on whether $\sigma$ is
real or complex, 
it follows from the commutative diagram following Definition~5.1 that
the volume of ${\rm vol}(\widetilde{\rm Pic}_F^0)$ is equal to  $2^{r_1}(2\pi\sqrt{2})^{r_2}/w_F$ times~${\rm Pic}^0_F$. This implies~{\sl (ii)}.

\beginsection 7. Reduced Arakelov divisors.

Let $F$ be a number field of degree~$n$.
In this section we introduce {\it reduced} Arakelov divisors associated to~$F$. These
form  a finite subset of ${\rm Div}_F^0$. The main result of this section is that the image of this set  in the groups  ${\rm Pic}^0_F$ and $\widetilde{\rm Pic}^0_F$ is in a certain sense regularly distributed. 

The results of this section extend work by Lenstra~[\LEX] and Buchmann and Williams~[\BUWINF] and
make certain statements by Buchmann in~[\BUJNTT,~\BUDPP,~\BUDU] more precise. In particular, Theorems~7.4 and~7.6
and Corollary~7.9
extend~ [\BUJNTT,~section~2], [\BUWINF,~Prop.2.7] and~[\BUDPP,~section~3.3]. 
Note  that 
in deducing the  corollaries below we did not make any particular effort to obtain the best possible estimates. They can most  certainly be improved upon.

Let $I$ be a fractional ideal. A non-zero element $f\in I$ is called {\it minimal} if 
 the only element $g\in I$ for which one has $|\sigma(g)|<|\sigma(f)|$ for all infinite
primes~$\sigma$, is~$g=0$. If $f\in I$ is minimal, then for every $h\in F^*$, the element $hf$ is minimal in the ideal~$h I$. In particular, if $h\in O_F^*$, the element $h f$ is minimal in the same ideal~$I$. Therefore there are, in general,  infinitely many minimal elements in~$I$.

If $D=(I,u)$ is an Arakelov divisor, then the minimal elements $f\in I$ are precisely the ones for which the 
open boxes $\{(y_{\sigma})_{\sigma}\in F_{\RR}: \hbox{$|y_{\sigma}|<|u_{\sigma}\sigma(f)|$ for all $\sigma$}\}$
contain only the point~$0$ of the lattice~$uI$. Note however that the notion of minimality depends only on~$I$ and is {\it
independent} of the metric induced by the element~$u$. {\it Shortest} elements $f\in I$
are the elements  for which $\abs{f}_D={\rm min}\{\abs{g}_D: \hbox{$g\in I-\{0\}$}\}$.  This notion depends on the
divisor $D=(I,u)$ and hence on the lattice~$uI$. It does not merely depend on~$I$.  Since $\abs{g}_D=\abs{ug}$ for each $g\in I$, the vectors $uf$ are the shortest non-zero vectors of the lattice
$uI$ associated to~$D$. The number of shortest elements in~$I$ is always finite.
Shortest vectors are  clearly minimal, but the converse is not true. It may even happen that 
a minimal element $f\in I$ is not a shortest element of the lattice $D=(I,u)$ for {\it any} choice of~$u$.
See section~9 for an explicit example.

\smallskip\noindent{\bf Definition.} An Arakelov divisor or oriented Arakelov divisor
$D$ in ${\rm Div}_F$ is called {\it reduced} if  it is of the form
$D=d(I)=(I,N(I)^{-1/n})$ for some fractional ideal~$I$, and if $1$ is a minimal element of~$I$. The set of
reduced Arakelov divisors is denoted by~${\rm Red}_F$. 

\medskip 
Since reduced Arakelov divisors have degree zero,
the covolume of the lattices associated to reduced Arakelov divisors  are equal to~$\sqrt{|\Delta_F|}$. With respect to the natural metric, 
$1\in O_F$ is a shortest and hence minimal element. Therefore the trivial
Arakelov divisor $(O_F,1)$ is  reduced. In general, if $D=d(I)$ is reduced,
the element $1\in I$ is merely
minimal and need not be a shortest element. However, the next proposition shows that 
it is not too far away from being so.

\proclaim Proposition 7.1.  Let $F$ be a number field of degree~$n$ and let $D=d(I)=(I,N(I)^{-1/n})$ 
be a reduced Arakelov divisor. Then we have
$$
\abs{1}_D\quad\le \quad\sqrt{n}\abs{x}_D,\qquad\hbox{for all non-zero $x\in I$.}
$$
In particular, the element $1\in I$ is at most $\sqrt{n}$ times as long as the shortest element in~$I$.

\smallskip\noindent{\bf Proof.} 
We have that $\abs{1}_D=\sqrt{n}N(I)^{-1/n}$. Since $1\in I$ is minimal,
every non-zero $x\in I$ has the property that $|\sigma(x)|\ge 1$ for some embedding $\sigma:F\longrightarrow\CC$. Therefore $\abs{x}_D\ge N(I)^{-1/n}|\sigma(x)|\ge N(I)^{-1/n}$. This proves the proposition.

\medskip
If $D=(I,u)$ is an Arakelov divisor and $f\in I$
is  minimal, then $1\in f^{-1}I$ is again minimal and the divisor $d(f^{-1}I)=(f^{-1}I,
N(fI^{-1})^{1/n})$ is reduced. In particular, if $f\in I$ is a shortest element,
then  the divisor $d(f^{-1}I)$ is reduced. Note however that
even though the element $1\in f^{-1}I$ is minimal,  it {\it  need not be} a shortest element.
Indeed,  even if it is true that $1$ is  a shortest vector of the lattice associated to
$(f^{-1}I,|f|^{-1}u)$,   it may not be a shortest vector of the lattice~$d(f^{-1}I)=
(f^{-1}I,N(fu^{-1})^{1/n})$,
which has a different metric. In
the section~9 we present an explicit example of this phenomenon. 

It is not so easy to say in terms of the associated ideal lattice $uI$ precisely what  it means
 that a divisor $D=(I,u)$ is reduced.
We make the following imprecise observation. When $1\in I$ is not merely minimal, but happens to be a shortest
element in~$I$,  then all roots of unity in~$F$ are  also shortest elements in~$I$.
Usually, these are the {\it only} shortest elements in~$I$. In that case the
arithmetic-geometric mean inequality implies that the Hermite constant
$\gamma(D)$, viewed as a function on ${\rm Pic}_F^0$ attains a local 
minimum at~$D=(I,N(I)^{-1/n})$. So, the lattices
corresponding to  reduced divisors  are the
``skewest" $O_F$-lattices around. But this holds only usually and locally.
\medskip
It is convenient to introduce the following notation.
\medskip\noindent{\bf Definition.} Let $F$ be a number field. Let $\Delta_F$ denote its discriminant and $r_2$ its number of complex infinite primes. Then we put
$$
\partial_F\quad=\quad \left({2\over{\pi}}\right)^{r_2}\sqrt{|\Delta_F|}.
$$

\proclaim Proposition 7.2. Let $F$ be a number field of degree~$n$. 
\item{(i)} Let $I$ be a fractional ideal. If
$d(I)=(I,N(I)^{-1/n})$ is a reduced Arakelov divisor,  then
the inverse $I^{-1}$ of $I$ is an $O_F$-ideal of norm at most~$\partial_F$.
\item{(ii)} The set ${\rm Red}_F$ of reduced Arakelov divisors is finite.
\item{(iii)}  The natural map ${\rm Red}_F\longrightarrow\widetilde{\rm Pic}^0_F$ is injective.

\smallskip\noindent{\bf Proof.} Since $1\in I$, the ideal $I^{-1}$ is contained in~$O_F$.
By Prop.~4.4~{\sl(ii)}  there exists a non-zero $f\in I$
for which $|N(I)^{-1/n}\sigma(f)|<\partial_F^{1/n}$
for each~$\sigma$.  Therefore, if
$N(I^{-1})>\partial_F$,  we have that
$|\sigma(f)|<1$ for each~$\sigma$, contradicting the minimality of~$1\in I$.
This proves~{\sl (i)}. Part~{\sl (ii)} follows at once from~{\sl (i)} and the fact that there are only finitely many
$O_F$-ideals  of bounded norm.

To prove {\sl (iii)}, suppose that the reduced Arakelov divisors $D=d(I)$ and
$D'=d(I')$ have the same image in $\widetilde{\rm Pic}^0_F$. Then there exists $f\in F^*$
so that $I'=fI$ and
$N(I')^{1/n}=N(I)^{1/n}f$. As in the proof of Prop.~5.3, it follows that all conjugates of~$f$ are
equal and hence that $f\in\QQ^*$. Since both $I$ and
$I'$ contain 1 as a minimal vector, this implies that
$f=\pm 1$. Since $f=N(I'I^{-1})^{1/n}>0$, we have that $f=1$ and hence $D=D'$ as required.

\medskip
Part~{\sl (iii)} of Proposition~7.2 does not hold when we replace
$\widetilde{\rm Pic}^0_F$ by~${\rm Pic}^0_F$.  See Example~9.3 below for
an example.  Incidentally, Theorem~7.7  below strengthens the statement considerably.

For every divisor $D=(I,u)$ of degree zero consider the following set of divisors of degree zero:
$$
\Sigma_D=\{(I,v)\in{\rm Div}^0_F:
\hbox{$\log(v_{\sigma})\le {1\over n}\log(\partial_F)$ for all $\sigma$}\}.
$$
If $N(I^{-1})\le \partial_F$, the set $\Sigma_D$ is a non-empty simplex. Indeed,
under this condition $\Sigma_D$ contains the divisor $(I,N(I)^{-1/n})$ and any  element of $\Sigma_D$ has the form
$(I,N(I)^{-1/n})+(O_F,w)$ with $w$ running over the exponentials of the vectors $y\in\left(\oplus_{\sigma}\RR\right)^0$
satisfying 
$$y_{\sigma}\le {1\over n}\left(\log(\partial_F)+\log(N(I))\right),\qquad
\hbox{for every~$\sigma$.}
$$ 
Since $\sum_{\sigma}{\rm deg}(\sigma)y_{\sigma}=0$, the set $\Sigma_D$ is a bounded simplex.

The following proposition expresses the notion of a reduced divisor in terms of these simplices.

\proclaim Proposition 7.3. An Arakelov  divisor of the form $D=d(I)=(I,N(I)^{-1/n})$ with $1\in I$, is reduced if and only if there is no other  divisor $D'=d(I')$ with $1\in I'$ for which the image of the simplex $\Sigma_D$
in~${\rm Pic}_F^0$ is contained in the image of~$\Sigma_{D'}$.

\smallskip\noindent{\bf Proof.} Suppose that $D=(I,N(I)^{-1/n})$ is reduced and that 
for some divisor $D'=(I',N(I')^{-1/n})$ with $1\in I'$ the image in~${\rm Pic}_F^0$ of $\Sigma_D$ is contained in the image of~$\Sigma_{D'}$.
Then $D$ and $D'$ lie on the same component of ${\rm Pic}_F^0$. This implies that  $I'=fI$ for some $f\in F^*$. Since $1\in I$ is minimal, so is $f\in {I'}$.
The simplex $\Sigma_{D'}$ or rather its image in ${\rm Pic}_F^0$ is equal to the set
$$\{(I,v/|f|):\hbox{$v$ satisfies $\log(v_{\sigma}/|\sigma(f)|)\le{1\over n}\log(\partial_F)$}\}.
$$ 
 Since $\Sigma_D\subset \Sigma_{D'}$, we have for each $\sigma$  that $\log(v_{\sigma}/|\sigma(f)|)\le 
{1\over n}\log(\partial_F)$ whenever $\log(v_{\sigma})\le {1\over n}\log(\partial_F)$. This implies that $|\sigma(f)|\ge 1$ for every~$\sigma$, contradicting the minimality of~$f\in I'$.

Conversely, suppose that $D=(I,N(I)^{-1/n})$ is not reduced. This means that $1\in I$ is not minimal.
Let $g\in I$ such that $|\sigma(g)|\le 1$ for all~$\sigma$. Consider the $O_F$-ideal $I'=g^{-1}I$. Then $\Sigma_D\subset \Sigma_{D'}$. Indeed, if 
$(I,v)\in \Sigma_D$ then $\log(v_{\sigma})\le {1\over n}\log(\partial_F)$ and hence  $\log(v_{\sigma}|\sigma(g)|)\le {1\over n}\log(\partial_F)$.
Since $(I,v)$ is equivalent to the divisor $(I',v|g|)$, this means precisely  $(I,v)$ is contained in~$\Sigma_{D'}$.

This proves the proposition.

\medskip
In the rest of this section we study the distribution of the image of the set
${\rm Red}_F$ in the compact groups ${\rm Pic}_F^0$ and~$\widetilde{\rm Pic}_F^0$
and estimate its size.
First we look at the image of the set ${\rm Red}_F$ in~${\rm Pic}_F^0$.  Theorem~7.4 
says that ${\rm Red}_F$ is rather {\it dense} in~${\rm Pic}_F^0$.  
\medskip

\proclaim Theorem 7.4. Let $F$ be a number field of degree~$n$ admitting $r_2$ complex infinite primes.  Then 
\item{(i)}  for any Arakelov  divisor $D=(I,u)$ of degree~0
there is a reduced divisor $D'$ and an element $f\in F^*$ so that 
$$D-D'=(f)+(O_F,v)$$ with
$$
\log|v_{\sigma}|\le {1\over n}\log(\partial_F).
$$
In particular, we have that
$$
\abs{D-D'}_{\rm
Pic}^{}\quod\le\log(\partial_F),\qquad\hbox{for each~$\sigma$;}
$$
\item{(ii)}
the natural map 
$$\mathop{\cup}\limits_{D}\Sigma_D\quad\longrightarrow\quad{\rm Pic}_F^0
$$ 
is surjective. Here $D$ runs over the reduced Arakelov divisors.

\smallskip\noindent{\bf Proof.} By Minkowski's Theorem (Prop.~4.4~{\sl (ii)}), there is a non-zero element
$f\in I$ satisfying $|u_{\sigma}\sigma(f)|\le\partial_F^{1/n}$ for every~$\sigma$.
Then there is also a shortest and hence a minimal such element~$f$. The divisor $D'=d(f^{-1}I)$ is then {\it reduced}.
It lies   on  the same component of ${\rm Pic}^0_F$ as~$D$. We have that
$$
D-D'=(f)+(O_F,v)
$$
where $v$ is the vector
$(v_{\sigma})_{\sigma}\in\prod_{\sigma}\RR^*_+$
with $v_{\sigma}=u_{\sigma}|\sigma(f)|N(f^{-1}I)^{1/n}$ and hence $\log|v_{\sigma}|=\log|u_{\sigma}\sigma(f)|+{1\over n}\log(N(f^{-1}I))$ for every~$\sigma$. Since $N(f^{-1}I)\le 1$,
this implies that $\log|v_{\sigma}|\le\log|u_{\sigma}\sigma(f)|$ which by assumption is at most
 $ {1\over n}\log(\partial_F)$ as required.

Since $\sum_{\sigma}\log(v_{\sigma})=0$,  Lemma~7.5 below implies that 
$$ \abs{D-D'}_{\rm Pic}^2=\abs{v}_{\rm Pic}^2\le n(n-1) \left({1\over n}\log(\partial_F)\right)^2.
$$
This proves~{\sl (i)}.
Part~{\sl (ii)} is merely a reformulation of part~{\sl (i)}. 
\medskip

\proclaim Lemma 7.5. Let $x_i\in\RR$ for $i=1,\ldots,n$. Suppose that $\sum_{i=1}^nx_i=0$ and that $x\in\RR$ has the property that $x_i\le x$ for all~$i=1,\ldots,n$. Then we have $\sum_{i=1}^nx_i^2\le n(n-1)x^2$.

\medskip
We 
leave the proof of this lemma to the reader.
The theorem says that ${\rm Pic}^0_F$ can be covered with   simplices $\Sigma_D$ centered in the reduced divisors~$D$.
We use  the Theorem to estimate the volume of the Arakelov class group~${\rm Pic}_F^0$ in terms of the number of reduced divisors.

\proclaim Corollary 7.6. Let $F$ be a number field of degree~$n$ with $r_1$ real and $r_2$ complex infinite primes. We have that
$$\eqalign{
 {\rm vol}({\rm Pic}_F^0)\quad&\le\quad
{{2^{-{{r_2}\over2}}n^{r_2-{1\over2}}}\over{(n-1)!}}\left(\log(\partial_F)\right)^{r_1+r_2} \#{\rm Red}_F,\cr
\quad&\le\quad(\log|\Delta_F|)^n\#{\rm Red}_F.\cr}
$$

\smallskip\noindent{\bf Proof.}  Let $D=d(I)=(I,N(I)^{-1/n})$ be reduced divisor. Then
the set $\Sigma_D$ is given by
$$
\Sigma_D=\{(I,N(I)^{-1/n})+(O_F,v):\log(v_{\sigma})\le  {1\over n}\left(\log(\partial_F)+\log(N(I))\right)\}.
$$
By Prop.~7.2~{\sl (i)} we have that $N(I^{-1})\le \partial_F$. 
This implies that the set $\Sigma_D$ is a non-empty simplex of volume  equal to
 $\left({1\over n}\log(\partial_F N(I))\right)^{r_1+r_2}$ times the volume of the standard simplex $\{(y_{\sigma})\in\mathop{\oplus}_{\sigma}\RR:\hbox{$\sum_{\sigma}y_{\sigma}=0$ and $y_{\sigma}\le 1$ for each~$\sigma$}\}$, which one  checks to be equal to~$2^{-{{r_2}\over2}}n^{n-{1\over2}}/(n-1)!$. This leads to the inequality
$$
 {\rm vol}({\rm Pic}_F^0)\quad\le\quad{{2^{-{{r_2}\over2}}n^{n-{1\over2}}}\over{(n-1)!}}\sum_D\left({1\over n}\log(\partial_FN(I))\right)^{r_1+r_2} $$
 Here the sum runs over the reduced divisors $D=(I,N(I)^{-1/n})$ of~$F$.

 Since $N(I)\le 1$, the first estimate follows follows. The second inequality follows by a rather crude 
 estimate from the first one. This proves the corollary.
 
\medskip
Next we prove some kind of converse to Theorem~7.4.
The following theorem and its corollary say that the image of the set ${\rm Red}_F$ is rather {\it
sparse} in the group~$\widetilde{\rm Pic}_F^0$. Recall that $F^*_+=\{x\in F^*:\hbox{$\sigma(x)>0$ for all real~$\sigma$}\}$.

\proclaim Theorem 7.7. Let $F$ be a number field. 
 \item{(i)} Let ${D}$ and ${D}'$ be
two reduced divisors in~$\widetilde{\rm Div}^0_F$. If 
there exists an element $f\in F^*_+$ for which
$$
{D}-{D}'+(f)=(O_F,v)
$$
with $|\log(v_{\sigma})|<{\rm log}({4\over3})$ for each $\sigma$, then  we have that
${D}={D}'$ in $\widetilde{\rm Div}^0_F$. 
Similarly, if $\abs{v}^{}_{\widetilde{\rm Pic}}<{\rm log}({4\over3})$,  we have that
${D}={D}'$ in $\widetilde{\rm Div}^0_F$. 
 \item{(ii)} The natural map 
$$\matrix{\mathop{\cup}\limits_{D'\in{\rm Red}_F}\{D'+(O_F,v):\,
\hbox{$v\in (F_{\RR,\,{\rm  conn}}^*)^0$ and $|\log(v_{\sigma})|< {1\over 2}{\rm log}({4\over 3})$
for each~$\sigma$}\}\cr\,\cr\Big\downarrow\cr\,\cr\widetilde{\rm Pic}_F^0\cr}
$$ 
is injective.

\smallskip\noindent{\bf Proof.}  Suppose that ${D}=d(I)$ and 
${D}'=d(I')$ are two reduced divisors with the property that
${D}-{D}'+(f)=(O_F,v)$ with $f\in F^*$ for which $\sigma(f)>0$ for all~$\sigma$. By Prop.5.3,  the images of $D$ and $D'$ in $\widetilde{\rm
Pic}_F^0$ lie on the same connected component of~$\widetilde{\rm Pic}_F^0$. We put $\lambda=N(I/I')^{1/n}$. Then  $\sigma(f)/\lambda=v_{\sigma}$.
Since $|{\rm log}(v_{\sigma})|<{\rm log}({4\over 3})$, we have that
$$
|{{\sigma(f)}\over{\lambda}}-1|=|v_{\sigma}-1|=|{\rm exp}({\rm log}(v_{\sigma}))-1|\le
{\rm exp}|{\rm log}(v_{\sigma})|-1<{\rm exp}({\rm log}(\hbox{${4\over 3}$}))-1=\hbox{${1\over 3}$},
$$
and hence
$$
|\sigma(f)-\lambda|\,\,\,<\,\,\,\hbox{${1\over3}$}\lambda,\qquad\hbox{for every~$\sigma$.}
$$ 
Since $D$ and $D'$ are reduced, the element 1 is
minimal in both
$I$ and~$I'$. Therefore both 1 and $f$ are minimal in $fI'=I$. 

If $\lambda$ is small, i.e., if $0<\lambda<{1\over2}$, we have that
$|\sigma(f)|\le|\sigma(f)-\lambda|+|\lambda|<{1\over3}\lambda+\lambda<{4\over3}\cdot{1\over2}<1$
for each~$\sigma$. In other words, $|\sigma(f)|<|\sigma(1)|$ for all~$\sigma$, contradicting  the
fact that $1\in I$ is minimal. If $\lambda$ is large, i.e., if $\lambda>{3\over2}$, we have  that
$|\sigma(f)|\ge|\lambda|-|\sigma(f)-\lambda|
\ge\lambda-{1\over 3}\lambda > {2\over 3}\cdot{3\over2}=1$ for each~$\sigma$. In other words,
$|\sigma(1)|<|\sigma(f)|$ for all~$\sigma$, contradicting the fact that $f\in I$ is a minimal
vector.

Therefore ${1\over2}\le\lambda\le{3\over2}$. This implies that $|\sigma(f-1)|\le
|\sigma(f)-\lambda|+|\lambda-1|<{1\over3}\lambda+|\lambda-1|\le{1\over 3}\cdot{3\over 2}+{1\over
2}=1=|\sigma(1)|$ for all~$\sigma$. Since  $1\in I$ is a minimal vector, this implies that
$f-1=0$. Therefore $I=I'$ and hence  $D=D'$. This proves the first statement.

If we know that $\abs{v}_{\widetilde{\rm Pic}}^{}<{\rm log}({4\over 3})$, then there is a 
totally positive unit $\varepsilon$ for which $|\log(\sigma(\varepsilon)v_{\sigma})|<{\rm log}({4\over 3})$
for each~$\sigma$.
Replacing $f$ by $\varepsilon f$ if necessary, we may then assume that 
$|\log(v_{\sigma})|<{\rm log}({4\over 3})$
for each~$\sigma$ and we are back in the earlier situation.
This proves~{\sl (i)}.

Part~{\sl(ii)} follows, because  {\sl (i)} implies that the sets 
$$\{D'+(O_F,v):\,\hbox{$v\in(F_{\RR,\,{\rm  conn}}^*)^0$ and $|\log(v_{\sigma})|< {1\over 2}{\rm log}({4\over 3})$
for each~$\sigma$}\}
$$ 
map injectively to $\widetilde{\rm Pic}_F^0$ and that their images are  mutually disjoint. This
proves the theorem.
\medskip

\proclaim Corollary 7.8. Let $F$ be a number field of degree~$n$.
Then we have that
$\#{\rm Red}_F\le {\rm vol}(\widetilde{\rm Pic}_F^0)\cdot (4n)^{{n\over 2}}$.

\smallskip\noindent{\bf Proof.}\ 
Theorem~7.7~{\sl (ii)} implies that the volume of $\widetilde{\rm Pic}_F^0$ is at least $\#{\rm
Red}_F$ times the volume of  the simplex $\{v\in((F_{\RR,\,{\rm  conn}}^*)^0:\hbox{$|\log(v_{\sigma})|< {1\over 2}{\rm log}({4\over 3})$ for each $\sigma$}\}$, which is equal to
$$
{{2^{-{{r_2}\over2}}n^{r_2-{1\over2}}}\over{(n-1)!}}\left(\hbox{$1\over 2$}\log(\hbox{$4\over 3$})\right)^{r_1+r_2}.
$$
Since this is at least $(4n)^{-{n\over2}}$, the result follows.

\medskip
\proclaim Corollary 7.9. There is a constant $c>0$, so that for every  number field $F$ of degree~$n$,
the number of reduced divisors
contained in a ball of radius~1 in~${\rm Pic}_F^0$ is at most $(cn)^{n/2}$.

\smallskip\noindent{\bf Proof.}  The reduced divisors whose images in ${\rm Pic}_F^0$ are contained in a ball of radius~1,
are contained in a subset  $S$ of $\widetilde{\rm Pic}_F^0$ of volume $2^{r_1}(2\pi\sqrt{2})^{r_2}/w_F$ times the volume of a unit ball  in ${\rm Pic}_F^0$.  By Theorem~7.7, the balls of radius $\hbox{${1\over 2}\log({4\over 3})$}$
centered in reduced divisors, are mutually disjoint in~$\widetilde{\rm Pic}_F^0$. Comparing the volume of the union of the disjoint balls with the volume of~$S$ leads to the estimate.

\medskip

\proclaim Corollary 7.10. Let $F$ be a number field of degree~$n$. 
Then we have that
$$
\left({\rm log}|\Delta_F|\right)^{-n}
\quad\le\quad{{\#{\rm Red}_F}\over{{\rm vol}({\rm Pic}^0_F)}} \quad\le\quad
\left(8\pi n\sqrt{2}\right)^{{n\over2}}.
$$

\smallskip\noindent{\bf Proof.}\   Since the volume of
$\widetilde{\rm Pic}^0_F$ is $2^{r_1}(2\pi\sqrt{2})^{r_2}/w_F$ times the volume of~${\rm Pic}_F^0$,
the inequalities follow from Corollaries~7.6 and~7.8 respectively. This proves the corollary
\medskip
We recall the following estimates for the volume of~${\rm Pic}_F^0$. They say that in a sense the
volume of ${\rm Pic}^0_F$ is approximately equal to~$\sqrt{|\Delta_F|}$.

\proclaim Proposition 7.11. Let $n\ge 1$. Then for every number field $F$  of degree~$n$ we have that
\item{(i)} 
$$
{\rm vol}({\rm Pic}^0_F)\le \sqrt{|\Delta_F|}\left({\rm log}{|\Delta_F|}\right)^{n-1};
$$
\item{(ii)} (GRH) there exists a constant $c>0$ only depending on the degree~$n$ so that
$$
{\rm vol}({\rm Pic}^0_F)\ge c\sqrt{|\Delta_F|}/{\rm log}{\rm log}{|\Delta_F|}.
$$

\smallskip\noindent{\bf Proof.} 
Part~{\sl (i)}  follows from Corollary~7.7, the fact that for every reduced divisor $d(I)$ the
ideal $I^{-1}$ is integral and has norm at most $\partial_F=\hbox{$\left({2\over{\pi}}\right)^{r_2}$}\sqrt{|\Delta_F|}$ and the estimate for the number of
$O_F$-ideals of bounded norm
provided by Lenstra in~[\LEB,~Thm.6.5]. 
Under assumption of the generalized Riemann Hypothesis (GRH) for the 
zeta function of the normal closure of~$F$, Buchmann and Williams~[\BUW,~(3.2)] obtained the estimate in~{\sl (ii)}. This proves the  proposition.

\beginsection 8. Quadratic fields.

Since the class group of $\QQ$ is trivial and since we have $\ZZ^*=\{\pm1\}$, the 
group ${\rm Pic}_{\QQ}^0$ is trivial and the degree map induces an isomorphism~${\rm Pic}_{\QQ}\cong\RR$. 
The narrow class group of~$\QQ$ is also trivial and it follows from Proposition~5.1 that
$\widetilde{\rm Pic}^0_{\QQ}=0$ and that $\widetilde{\rm Pic}_{\QQ}^{}$ is isomorphic to~$\RR^*$.

This is the whole story as far as $\QQ$ is concerned. In this section we briefly work out the theory of the previous sections for quadratic number
fields. For these fields the language of binary quadratic forms is often used~[\LEX,~\SH].

\medskip\noindent
{\bf Example 8.1.} 
For complex quadratic fields~$F$, the torus $T^0$ of section~2 is trivial so that 
the group ${\rm Pic}^0_F$ is
canonically isomorphic to the class group $Cl_F$ of~$F$. The group $\widetilde{\rm Pic}^0_F$ is an
extension  of $Cl_F$ by a circle group of length~$2\pi\sqrt{2}/w_F$. Here $w_F=2$ except when $F=\QQ(i)$ or $\QQ({{1+\sqrt{-3}}\over2})$, in which case $w_F=4$ or~6 respectively.

We describe the reduced
Arakelov divisors of~$F$. 
Let $D=(I,N(I)^{-1/2})$ be reduced. The fact that $1$ is a minimal element of~$I$
simply means that it is a shortest vector in the corresponding lattice 
in~$F_{\RR}\cong\CC$.  We write
$I=\ZZ+f\ZZ$ for some $f$ in the upper half plane~$\{z\in\CC:{\rm Im}(z)>0\}$. 
Since $O_F\cdot I\subset I$, we have that
$f={{b+\sqrt{\Delta_F}}\over{2a}}$ for certain $a,b\in\ZZ$, $a>0$ and 
that $b^2-4ac=\Delta_F$ for some~$c\in\ZZ$. The $O_F$-ideal $I^{-1}$ is generated by
$a$ and~${{b+\sqrt{\Delta_F}}\over{2}}$ and has norm~$a$. 
For complex quadratic fields, the simplices  $\Sigma_D$ introduced in section~6 are simply points.

Since $f$ is unique up to addition of  an integer, the ${\rm SL}_2(\ZZ)$-equivalence class of the
binary quadratic form $N(X+fY)/N(I)=aX^2+bXY+cY^2$ is well defined. The form has
discriminant~$\Delta_F$. If we choose $f$ to lie in the usual fundamental domain
for the action of ${\rm SL}_2(\ZZ)$ on the upper half plane, the corresponding quadratic form
is reduced in the sense of Gauss. There is a slight ambiguity here. If $|f|=1$, the reduced Arakelov
divisors $d(\ZZ+f\ZZ)$ and $d(\ZZ+\overline{f}\ZZ)$ give rise to
the quadratic forms $aX^2+bXY+aY^2$ and $aX^2-bXY+aY^2$ respectively. If $f$ is not a root of unity, the
Arakelov divisors are distinct, but the two quadratic forms are ${\rm
SL}_2(\ZZ)$-equivalent and only one of them is reduced.  Apart from this ambiguity, the map that
associates to a reduced Arakelov divisor its associated reduced quadratic form, is a bijection.

\medskip
\noindent{\bf Example 8.2.}
Any real quadratic field $F$ can be written as $\QQ(\sqrt{\Delta_F})$ where $\Delta_F$ denotes the discriminant
of~$F$. The group
${\rm Pic}^0_F$ is an extension of the class group by a circle group and the group $\widetilde{\rm
Pic}^0_F$ is an extension of the {\it narrow} class group by a circle group. We describe the reduced Arakelov
divisors of
$F$. Let $\sigma$ and $\sigma'$ denote the two infinite primes of~$F$. To be definite, we let $\sigma$
denote the embedding that maps $\sqrt{\Delta_F}$ to the positive square root of $\Delta_F$
in~$\RR$.
Let $D=d(I)=(I,N(I)^{-1/2})$ be reduced. The fact that $1\in I$
is minimal implies that we can write $I=\ZZ+f\ZZ$ for a
unique $f$ satisfying $\sigma(f)>1$ and $-1<\sigma'(f)< 0$. The fact that we have
$O_F\cdot I\subset I$ implies that $f={{b+\sqrt{\Delta_F}}\over{2a}}$ where $\Delta_F=b^2-4ac$ for
some $c\in\ZZ$. The conditions on~$\sigma(f)$ and $\sigma'(f)$ say that $a>0$ and
$|\sqrt{\Delta_F}-2a|<b<\sqrt{\Delta_F}$.  The $O_F$-ideal $I^{-1}$ is generated by $a$
and~${{b+\sqrt{\Delta_F}}\over{2}}$. Its norm is~$a$. The simplex $\Sigma_D$ of section~6
is an interval of length $\sqrt{2}\,\log({{\sqrt{\Delta_F}}\over a})$ centered in $D$.

The map that
associates the quadratic form
$aX^2+bXY+cY^2$ to the reduced divisor $D=(I,N(I)^{-1/2})$, is a bijection between the set of
reduced Arakelov divisors of $F$ and the set of reduced binary quadratic forms of
discriminant~$\Delta_F$ with $a>0$.

The element $1\in I$ 
is a shortest vector
precisely when both $\abs{f}$ and $\abs{f-1}$ are at least~$\abs{1}=\sqrt{2}$. This condition is not
always satisfied. Drawing a picture, one sees that it is when $\sigma(f)-\sigma'(f)\ge 2$, or equivalently when
$a<{1\over 2}\sqrt{\Delta_F}$, but this is not a necessary condition. 

When $D=d(I)$ and  $I=\ZZ+f\ZZ$ as above, then the  vector $f$ is a minimal element
of~$I$. Therefore $D'=d(f^{-1}I)$ is a reduced Arakelov divisor. We have that
$D=D'+(f)+(O_F,v)$, where $v\in F^*_{\RR}\cong\RR^*\times\RR^*$ is the vector
$(|\sigma'(f)/\sigma(f)|^{1/2}, -|\sigma(f)/\sigma'(f)|^{1/2})$. The distance between the images
of $D$ and
$D'$ in ${\rm Pic}_F^0$ is equal to~$\abs{v}_{\rm Pic}$. Since $f={{b+\sqrt{\Delta}}\over{2a}}$, we have
$\abs{v}_{\rm Pic}={1\over {2\sqrt{2}}}{\rm log}|{{b+\sqrt{\Delta_F}}\over{b-\sqrt{\Delta_F}}}|$. In this way we recover Lenstra's
distance formula~[\LEX,~(11.1)]. The divisor $D'$ is the `successor' of $D$ in its component, in the sense
that there are no reduced divisors on the circle between $D$ and $D'$. In order to obtain $D$'s `predecessor', take $g$
the shortest minimum  such that
$|\sigma(g)|<|\sigma'(g)|$.  Then the Arakelov divisor $d({g}^{-1}I)$ is the predecessor
of~$D$.

Lenstra's group $\cal F$ or rather its topological completion $\overline{\cal F}$, is closely related to the oriented
Arakelov class group of the real quadratic field~$F$ and several of the results in~[\LEX] 
are special cases of  the results
of this paper.  The group $\overline{\cal F}$ is not
quite {\it equal} to~$\widetilde{\rm Pic}_F^0$ but it admits a degree~2 cover onto it. More
generally, for a number field $F$ we let
${\rm Pic}_{F}^+$ denote  the group $\widetilde{\rm Div}_F^0$ modulo its subgroup~$\pm F^*_+$. When
$F$ is totally complex, i.e., when $r_1=0$, this is simply~$\widetilde{\rm Pic}_F^0$. When $r_1>0$
however, there is an exact sequence
$$
0\longrightarrow\quad\{\pm 1\}^{r_1}/\{\pm 1\}\quad\longrightarrow\quad {\rm
Pic}_{F}^+\quad\longrightarrow\quad\widetilde{\rm Pic}_F^0 \quad\longrightarrow 0.
$$
Let $(F^*_{\RR})^0=\{u\in F_{\RR}:|N(u)|=1\}$. The topological structure of ${\rm Pic}_{F}^+$ can be seen from the 
exact sequence
$$
0\longrightarrow\quad(F^*_{\RR})^0/\pm
O_{F,+}^*\quad\longrightarrow\quad{\rm Pic}_{F}^+\quad\longrightarrow\quad
Cl_{F,+}\quad\longrightarrow 0,
$$
realizing ${\rm Pic}_{F}^+$ as an extension of the narrow class group $Cl_{F,+}$ by a
$2^{r_1-1}$-component  Lie group. When $F$ is real quadratic, the group ${\rm Pic}_{F}^+$ is equal to
Lenstra's group $\overline{\cal F}$.

\beginsection 9. Reduced Arakelov divisors; examples and counterexamples.

Let $F$ be a number field of degree~$n$ and discriminant~$\Delta_F$.
Theorems~7.4 and~7.7 say
that the image of the set ${\rm Red}_F$ of reduced Arakelov divisors is, in a precise sense, 
rather regularly distributed in the groups ${\rm Pic}^0_F$ and $\widetilde{\rm Pic}^0_F$. 
In this section we discuss these results and we consider variations in the 
definition of the set of reduced divisors.

Theorem~7.4 says that the image of ${\rm Red}_F$ is rather `dense' in~${\rm Pic}_F^0$.
I do not know how to prove similar result for the larger group
$\widetilde{\rm Pic}_F^0$ rather than ${\rm Pic}_F^0$.  I
cannot even exclude the possibility that some components of~$\widetilde{\rm Pic}_F^0$ contain no
reduced Arakelov divisors at all. The problem is related to the following question.

\medskip\noindent{\bf Question 9.1.}  Let
$L\subset\RR^n$ be a lattice of covolume~1. Suppose that all non-zero vectors of $L$ have all
their coordinates different from zero and have length at least~$\varepsilon$. Does there exist a
minimal vector $(x_i)\in L$ with
$x_i>0$ for all~$i$? Here a vector $(x_i)$ is called minimal if the only vector $(y_i)\in\RR^n$
with $|y_i|<|x_i|$ for all~$i$ is the zero-vector.
Does there exist such a vector of length
$O(\varepsilon^{-N})$ for some $N$ that only depends on~$n$? Can one compute it efficiently?

\medskip\noindent
In the other direction, Theorem~7.5 implies that the image of ${\rm Red}_F$ in
$\widetilde{\rm Pic}_F^0$ is rather `sparse'. When we replace $\widetilde{\rm Pic}_F^0$
by ${\rm Pic}_F^0$, the theorem is no longer true.
First of all the map ${\rm Red}_F\longrightarrow{\rm Pic}_F$ is in general not injective.
In addition, it may happen that distinct reduced divisors have images in ${\rm Pic}_F^0$ that
are  much closer to one another than the bound ${\rm log}({4\over3})$ of Theorem~7.5. However,
by Corollary~7.9, the number of reduced divisors in a ball in~${\rm
Pic}_F^0$ of radius~$1$ is bounded by a constant only depending on  the degree of~$F$.

\proclaim Lemma 9.2. Let $F$ be a number field of degree~$n$, let $D=(I,u)$ be an
Arakelov divisor and suppose~$f\in I$. Then
\item{(i)} $d(f^{-1}I)=d(I)$ in ${\rm Div}_F^0$ if and only if $f$ is a unit of~$O_F$.
\item{(ii)} The classes of $d(f^{-1}I)$ and $d(I)$ in ${\rm Pic}_F^0$ are equal if and only if
$f$ is the product of a unit and an element $g\in F^*$  all of whose absolute values
$|\sigma(g)|$ are equal.
\item{(iii)}  $\abs{d(I)-d(f^{-1}I)}_{\rm Pic}^{}<2\sqrt{n}\,{\rm max}_{\sigma}|{\rm log}|\sigma(f)||$.

\smallskip\noindent{\bf Proof.}\  Part~{\sl (i)} follows from the fact that $I=f^{-1}I$
if and only if we have $f\in O_F^*$. Since we have
$$
d(f^{-1}I)-d(I)=(fO_F,|N(f)|^{1/n}),
$$
the class of this  divisor is trivial in~${\rm Pic}_F^0$ if and only if
there is $g\in F^*$ for which $f=\varepsilon g$ for some unit $\varepsilon\in O_F^*$
and $|\sigma(g)|^{-1}=|N(f)|^{-1/n}$ for all~$\sigma$. Since $|N(g)|=|N(f)|$,
the second relation is equivalent to the fact that the $|\sigma(g)|$ are all equal. This proves~{\sl (ii)}.

To prove~{\sl (iii)} we note that $\abs{d(I)-d(f^{-1}I)}_{\rm
Pic}^{}\le\sqrt{n}\, {\rm max}_{\sigma}|{\rm log}|\sigma(f)/N(f)^{1/n}||$ which is at most 
$\sqrt{n}$ times
${\rm max}_{\sigma}|{\rm log}|\sigma(f)||+{1\over n}\sum_{\sigma}{\rm deg}(\sigma)|{\rm
log}|\sigma(f)|$. This easily implies the estimate.

 This implies that
$f$ is a root of unity. This completes the proof of the lemma.

\medskip
Prop.~7.2~{\sl (iii)} says that the natural map from the set of reduced divisors ${\rm Red}_F$ to the oriented Arakelov class group~$\widetilde{\rm Pic}_F^0$
is injective. The following example shows that, in general, the map 
${\rm Red}_F\longrightarrow{\rm Pic}_F^0$ is {\it not}.
\medskip
\noindent{\bf Example 9.3.} 
Let $a>b\ge 1$ and put $\Delta=b^2-4a^2$. Suppose that $\Delta$ is
squarefree and let $F$ denote the complex quadratic number field 
$\QQ(\sqrt{\Delta})$. Let $I$ denote the fractional 
$O_F$-ideal $\ZZ+f\ZZ$ where $f={{b+\sqrt{\Delta}}\over {2a}}$. Then $1\in I$ is minimal. Let
$\sigma:F\longrightarrow\CC$ denote the unique infinite prime. Since
$\sigma(f)$ has absolute value~1, the element $f$ is also minimal. 
Since $f$ is not a unit of~$O_F$, Lemma~9.2 implies that the reduced divisors $d(I)$ and
$d(f^{-1}I)$ are distinct, but that their classes in~${\rm Pic}_F^0$ are equal.

\medskip
Theorem~7.7 says that the distance between the images of the reduced divisors in~$\widetilde{\rm Pic}_F^0$
is bounded from below by an absolute constant. The following example shows that this is false for the
Arakelov class group~${\rm Pic}_F^0$.
\medskip
\noindent{\bf Example 9.4.}  Let $n$ be a large even integer such that
$\Delta=n^2+1$ is squarefree and consider the field $F=\QQ(\sqrt{\Delta})$. 
Let $f={{1+\sqrt{\Delta}}\over n}\in F$. Then 1 is a minimal element in
$I=\ZZ+f\ZZ$. The conjugates $\sigma(f)$ are close to~1 and~$-1$ respectively. Indeed, we have 
$|{\rm log}|\sigma(f)||\approx{\Delta}^{-1/2}$ for each infinite prime~$\sigma$. It follows
from Lemma~9.2~{\sl (iii)} that the classes of the reduced divisors $d(I)$ and $d(f^{-1}I)$ are at distance at most
$2\sqrt{2}\,{\Delta}^{-1/2}$ in~${\rm Pic}_F^0$.
\medskip

The definition of the set ${\rm Red}_F$ is rather delicate
as we'll see now by considering slight variations of it. We let
${\rm Red}'_F$ denote the set of divisors $d(I)$ for which $1\in I$ is a  {\it shortest}
rather than a {\it minimal}  vector and write ${\rm Red}''_F$
for the set of divisors $d(I)$ for which we have 
$N(I^{-1})\le\partial_F=\hbox{$\left({2\over{\pi}}\right)^{r_2}$}\sqrt{|\Delta_F|}$
and  for which $1\in I$ is merely {\it
primitive}, i.e., not divisible by an integer $d\ge 2$.  
Since {\it shortest} implies {\it minimal} and {\it minimal} implies {\it primitive},
we have the following inclusions of finite sets
$$
{\rm Red}'_F\quad\subset\quad{\rm Red}_F\quad\subset\quad{\rm Red}''_F.
$$
Theorem~7.4 says that the set  ${\rm Red}_F$ is rather `dense' in the Arakelov divisor
class group. 
It is not clear whether the set ${\rm Red}'_F$ has the same property. 
The  proof of Theorem~7.4 showing that every $D=(I,u)$ of degree~0 is close to a
reduced divisor~$D'\in {\rm Red}_F$ does not work for~${\rm Red}'_F$. Indeed, tracing the steps of the proof
of Theorem~7.4, we see that if  $f\in I$ is a shortest vector, it is also minimal and hence the element $1\in f^{-1}I$ is
{\it minimal}. It follows that the divisor $d(f^{-1}I)$ is in~${\rm Red}_F$. However, $1$ need not be a {\it
shortest} vector in $f^{-1}I$ so that
$d(f^{-1}I)$ may not be contained in~${\rm Red}'_F$. 

\medskip
The following example shows that this phenomenon actually occurs. It shows that  the set ${\rm Red}'_F$ 
is, at least in this  sense, too small.

\medskip\noindent{\bf Example 9.5.} We present examples of  reduced Arakelov divisors $D=d(I)$ with
the property that the element $1\in I$ is {\it not} a shortest vector of
the lattice $I$ associated to $(I,u)$ for {\it any} $u\in F_{\RR}^*$. 
This implies that $D$ is not equal to $d(f^{-1}J)$ for
any divisor $D'=(J,v)$ and a shortest element $f\in J$. Indeed, if that were the case,
$1$ would be shortest vector in the lattice associated to the Arakelov divisor $(I,f^{-1}v)$.

Let $F$ be a real quadratic number field of discriminant~$\Delta$.
Then $F=\QQ(\sqrt{\Delta})$. Suppose that  $d(I)$ is a reduced Arakelov divisor. We write $I=\ZZ+f\ZZ$ where
$f>0$ and $-1<\overline{f}<0$. Here we identify $F$ with its image in~$\RR$ through one of its
embeddings and we write $f\mapsto\overline{f}$ for the other embedding. 

\smallskip\noindent{\bf Claim.} If $N(f-{1\over2})>-{3\over 4}$, then $1$ is not a shortest element of~$I$ for {\it any} degree zero Arakelov divisor~$(I,u)$.
\smallskip
\noindent{\bf Proof.}  Suppose that $D=(I,u)$ has degree~0. Then we have
$u=(\upsilon/\sqrt{N(I)},\upsilon^{-1}/\sqrt{N(I)})$ for some $\upsilon\in\RR^*_{>0}$. 
Suppose that $1\in I$ is a shortest vector
in the lattice associated to~$D$.
This implies in particular that
$\abs{1}_D^{}\le\abs{f}_D^{}$ and $\abs{1}_D^{}\le\abs{f-1}_D^{}$. This means that
$\upsilon^{-2}+\upsilon^2\le \upsilon^{-2}f^2+\upsilon^2\overline{f}^2$
and that $\upsilon^{-2}+\upsilon^2\le \upsilon^{-2}(f-1)^2+\upsilon^2(\overline{f}-1)^2$.  In other words we have that
$\upsilon^4\le (f^2-1)/(1-\overline{f}^2)$ and $\upsilon^4\ge (2f-f^2)/(\overline{f}^2-2\overline{f})$ respectively.
Therefore, if the upper bound for $\upsilon^4$ is smaller than the lower bound, there cannot exist such~$\upsilon$.
This happens precisely when $(f-\overline{f})(2f\overline{f}-f-\overline{f}+2)> 0$. Since
$f-\overline{f}$ is positive, this means that $2f\overline{f}-f-\overline{f}+2>0$ which is equivalent to
$N(f-{1\over2})>-{3\over 4}$. This proves the claim.

When
$f={{b+\sqrt{\Delta}}\over{2a}}$ as in section~8, a sufficient condition for the inequality of the claim 
to hold is that~$a\ge{1\over{\sqrt{3}}}\sqrt{\Delta}$.
An explicit example is provided by the field $\QQ(\sqrt{21})$ and the 
reduced divisor $d(I)$ where $I=\ZZ+f\ZZ$ with
$f={{3+\sqrt{21}}\over6}$.

\smallskip

In the other direction, it may happen that the image of  ${\rm Red}''_F$ is very dense
in~$\widetilde{\rm Pic}_F^0$, so that an analogue of Theorem~7.7 does not hold for this set.  We
present two examples, due to H.W.~Lenstra showing that for some 
number fields certain small open balls in
$\widetilde{\rm Pic}^0_F$ contain the images of very many $D\in{\rm Red}''_F$. Both examples exploit the
existence of certain `very small' elements in~$F$. In the first example these are contained in a proper
subfield, but this is not the case in the second example.

\medskip\noindent{\bf Example 9.6.}
Let $F$ be a number field of degree~$n$ containing~$\QQ(i)$. Let $m,\, m'\in\ZZ$ satisfy ${1\over
2}|\Delta_F|^{1/2n}<m,\, m'<|\Delta_F|^{1/2n}-1$. Let
$I$ and $I'$ denote the inverses of the $O_F$-ideals generated by
$m-i$ and~$m'-i$ respectively. Then $1$ is primitive in both $I$ and~$I'$ and the norms of $I^{-1}$
and ${I'}^{-1}$ do not exceed $\partial_F=\left({2\over{\pi}}\right)^{r_2}|\Delta_F|^{1/2}$. 
It follows that $d(I)$ and $d(I')$ are in~${\rm
Red}''_F$. If the images of $d(I)$ and
$d(I')$ in $\widetilde{\rm Pic}^0_F$ are equal, Prop.~6.4 implies that $I=mI'$ for some $m\in\QQ^*$.
Since $1$ is primitive in both $I$ and $I'$, it follows that~$m=\pm 1$. This implies that
$I=I'$ and hence that $N(I)=m^2+1$ is equal to $N(I')={m'}^2+1$, so that~$m=m'$.
Therefore $d(I)$ and $d(I')$ are distinct in $\widetilde{\rm Pic}^0_F$, whenever $m$ and $m'$ are.

Assume in addition that  $|m-m'|< |\Delta_F|^{1/3n}$ and that $|\Delta_F|>4^{6n}$. Then
the distance between $m$ and $m'$ is much smaller than $m$ and~$m'$ themselves.
The  distance between the Arakelov divisors $d(I)$ and $d(I')$ in $\widetilde{\rm Pic}^0_F$
is at most $\sqrt{n}\,|{\rm log}((m-i)/(m'-i))|$. This does not exceed $\sqrt{n}|m-m'|/(|m-i|-|m-m'|)\le 
\sqrt{n}||\Delta_F|^{1/3n}/({1\over2}|\Delta_F|^{1/2n}-|\Delta_F|^{1/3n})\le 4\sqrt{n}||\Delta_F|^{-1/6n}$.

In this way we obtain $|\Delta_F|^{1/3n}$ elements of ${\rm Red}''_F$ whose images in  
$\widetilde{\rm Pic}^0_F$ are distinct, but are as close as $4\sqrt{n}||\Delta_F|^{-1/6n}$ to
one another. By varying $F$ over degree $n/2$ extensions of $\QQ(i)$, we can make $|\Delta_F|$ as large
as we like. One may replace $\QQ(i)$ by any number field and proceed similarly.

\medskip\noindent{\bf Example 9.7.}
Let $n\ge 4$ and $a\in\ZZ$  be such that the polynomial $X^n-a$ is irreducible
over~$\QQ$. Let $\alpha$ denote a zero and put $F=\QQ(\alpha)$. Suppose that the ring of integers
of $F$ is equal to $\ZZ[\alpha]$. There are infinitely many such integers~$a$. Then
$|\Delta_F|=n^n|a|^{n-1}$ and $|\sigma(\alpha)|=|a|^{1/n}$ for every infinite prime~$\sigma$.
Let $m,m'\in\ZZ $ satisfy ${1\over 2}|a|^{1/2-1/2n}+|a|^{1/n}<m,m'<|a|^{1/2-1/2n}$ and
$|m-m'|\le|a|^{1/4}$. Consider two  Arakelov divisors $d(I)$ and
$d(J)$ given by $I^{-1}=(m-\alpha)O_F$ and ${J}^{-1}=(m'-\alpha)O_F$. The norms of $I^{-1}$ and
$J^{-1}$ are at most $\partial_F=\left({2\over{\pi}}\right)^{r_2}|\Delta_F|^{1/2}$.  
Since both  $I$ and $J$
contain $1$ as a primitive element,  we have $d(I), d(J)\in{\rm Red}''_F$. The argument used in Example~9.6 shows that
the images of $d(I)$ and
$d(J)$ in $\widetilde{\rm Pic}^0_F$ are distinct when  $m\not=m'$. The difference between
$d(I)$ and $d(J)$  is equal to
$(I{J}^{-1},N(I{J}^{-1})^{1/n})$ which is equivalent to $(O_F,v)$ where
$$
v={{m-\sigma(\alpha)}\over{m'-\sigma(\alpha)}}\left|
N\left({{m'-\alpha}\over{m-\alpha}}\right)\right|^{1/n}.
$$
It follows that $\abs{d(I)-d(J)}_{\widetilde{\rm Pic}}$ is at most $2\sqrt{n}\,{\rm max}_{\sigma}\left|{\rm
log}({{m-\sigma(\alpha)}\over{m'-\sigma(\alpha)}})\right|$.
Since $(m-\sigma(\alpha))/(m'-\sigma(\alpha))=1+(m-m')/(m'-\sigma(\alpha))$ 
and since
$|m'-\sigma(\alpha)|\ge m'-|\sigma(\alpha)|\ge{1\over
2}|a|^{1/2-1/2n}$, the absolute value of the logarithm of
$(m-\sigma(\alpha))/(m'-\sigma(\alpha))$ is at most
$4|m-m'|/|a|^{1/2-1/2n}$ for each~$\sigma$. It follows that $\abs{d(I)-d(J)}_{\widetilde{\rm Pic}}^{}$
is at most
$4 \sqrt{n}|a|^{-1/4+1/2n}$, which becomes arbitrarily small as $|a|$ grows.

\beginsection 10. Computations with reduced Arakelov divisors.

In this section we discuss the set of reduced Arakelov 
divisors from a computational point of view. Our presentation is rather informal.
In particular, we do not say  much about the  accuracy of the 
approximations required to perform the computations  with the real and complex numbers
involved.
See~[\THIE] for a more rigorous approach.
Since Arakelov divisors can be represented as
lattices in the Euclidean space $F_{\RR}$, lattice reduction algorithms play an
important role. When the degree $n$ of the number field is large, the celebrated 
LLL-reduction algorithm~[\LLL, \LENSH] is an important  tool.

We suppose that the number field $F$ is given as
$\QQ(\alpha)$ where $\alpha$ is the zero of some irreducible monic polynomial $\varphi(X)\in\ZZ[X]$.
We assume that we have
already computed an LLL-reduced 
basis $\{\omega_1,\ldots,\omega_n\}$ for the ring of integers~$O_F$ embedded in~$F_{\RR}$. In
other words, we have an explicit lattice
$$
O_F=\omega_1\ZZ+\ldots+\omega_n\ZZ\subset F_{\RR},
$$
with, say, an LLL-reduced basis $\{\omega_1,\ldots,\omega_n\}$.
Such a basis can be computed as explained in~[\LEB, sect.4] or~[\CO,~sect.6.1] combined with a
basis reduction algorithm.  We have also computed a multiplication table i.e.,
coefficients $\lambda_{ijk}\in\ZZ$ for which $\omega_i\omega_j=\sum_k\lambda_{ijk}\omega_k$. 
The discriminant $\Delta_F$ of~$F$ is the integer given by $\Delta_F={\rm det}({\rm
Tr}(\omega_i\omega_j))$.    By~~[\LEB,~section 2.10] we have that $\lambda_{ijk}= |\Delta_F|^{O(n)}$. 
We view the degree $n$ of $F$ as fixed and estimate the running times
of the algorithms in terms of $|\Delta_F|$. 

An Arakelov divisor or oriented Arakelov divisor
$D=(I,u)$ is determined by its associated ideal~$I$ and the
vector~$u\in F^*_{\RR}\cong\prod_{\sigma}F_{\sigma}^*$. It can be represented by  an
$n\times n$ matrix $\lambda_{ij}$ for which the vectors
$\sum_{ij}\lambda_{ij}\omega_j$ are an LLL-reduced basis for the lattice
$I\subset F_{\RR}$, together with a sufficiently accurate
approximation of the vector $u=(u_{\sigma})_{\sigma}$.  We have that $\lambda_{ij}=O(N(I))$. See~[\THIE]. 
In practice, one might want to take logarithms and work
with the vectors~$({\rm log}(u_{\sigma}))_{\sigma}$. There are efficient algorithms to multiply ideals, to
compute inverses and to test for equality. See~[\CO,~sects.4.6-8]. Using these one can compute
efficiently in the group ${\rm Div}_F$. The algorithms have been 
implemented in the LiDIA, MAGMA and PARI software packages~[\LID, \MAG, \PARI]. 

Rather than the Arakelov divisor group, we are  interested in computing in the {\it Arakelov 
class group}~${\rm Pic}_F^0$. We do calculations in this group by means of the set
${\rm Red}_F$ of {\it reduced} divisors
in~${\rm Div}^0_F$. By Theorems~7.4 and~7.7, the image of the finite set ${\rm Red}_F$ is in
a certain sense regularly distributed in the compact groups ${\rm Pic}_F^0$ and~$\widetilde{\rm Pic}_F^0$. 
Reduced divisors have one further property that is important for our application: a reduced divisor $D$
 is of the form $D=d(I)=(I,N(I)^{-1/n})$ where
$I^{-1}$ is and integral ideal of norm at most~$\partial_F=\left({2\over{\pi}}\right)^{r_2}|\Delta_F|^{1/2}$. 
Therefore $D$ can be represented using
only~$({\rm log}|\Delta_F|)^{O(n)}$ bits.  

Before describing the algorithms, we formulate a lemma concerning the Lenstra-Lenstra-Lovasz (LLL)
lattice reduction algorithm~[\LLL].

\proclaim Lemma 10.1. Let 
$\b_1,\ldots,\b_n$  be an  LLL-reduced basis of a real vector space~$V$. Then for every vector ${\bf x}=\sum_{i=1}^nm_i\b_i$  of~$V$ we have that
$$
|m_i|\abs{{\b}^*_i}\quad\le\quad\left({3\over {\sqrt{2}}}\right)^{n-i}\abs{{\bf x}},\qquad\hbox{for
$1\le i\le n$.}
$$
Here $\b_1^*,\ldots, \b_n^*$ denotes the  Gram-Schmidt orthogonalization of the basis~$\b_1,\ldots,\b_n$.

\smallskip\noindent{\bf Proof.} See~[\LENSH].

\proclaim Corollary 10.2. 
Let $\b_1,\ldots,\b_n$ be an LLL-reduced basis of a real vector space~$V$.
Then we have for any  vector ${\bf x}=\sum_{i=1}^nm_i\b_i$ in $V$ that
$$
|m_{i}|\quad\le \quad 2^{{n-1}\over 2}\left({3\over 2}\right)^{n-i}{{\abs{\bf x}}\over{\abs{\b_1}}},\qquad
\hbox{for $i=1,\ldots,n$.}
$$

\smallskip\noindent
{\bf Proof.} The LLL-conditions~[\LLL]  imply  $\abs{\b_1^*}\le 2^{{i-1}\over2}\abs{\b_i^*}$ for  every $i=1,2\ldots,n$. 
 Since $\b_1=\b_1^*$, the result  follows from Lemma~10.1. 

\medskip
We have the following basic algorithms at our disposal.
For number fields of fixed degree~$n$, each runs in time polynomial in~${\rm log}|\Delta_F|$.

\proclaim Algorithm 10.3. ({\it Reduction algorithm}) Given an Arakelov divisor $D=(I,u)\in{\rm Div}_F^0$, 
\smallskip
\item{--} check whether it is reduced  or not;
\item{--} compute a reduced divisor $D'$ that is close to~$D$ in~${\rm Pic}_F^0$. 

\smallskip\noindent{\bf Description.} We compute an LLL-reduced basis $\b_1,\ldots,\b_n$ of the
lattice ~$uI\subset F_{\RR}$. Any  shortest vector  ${\bf x}=\sum_{i=1}^nm_i\b_i$ in the lattice
satisfies $\abs{{\bf x}}/\abs{{\bf b}_1}\le 1$. Therefore
Corollary~10.2 implies that the coordinates $m_i\in\ZZ$ are bounded independent of the
discriminant of~$F$. To compute a shortest vector in the lattice in
time polynomial in~$\log|\Delta_F|$, we may  therefore just try all possible~$m_i$. 

In order to compute a  reduced divisor $D'$ that is close to~$D$ in~${\rm Pic}_F^0$, 
we compute   a shortest vector $f$ in the lattice
$I$ associated to~$D$. The  divisor $D'=d(f^{-1}I)$ is then reduced. Moreover, by Theorem~7.4 or rather its proof,
the divisor $D'$ has the property that
$\abs{D-D'}_{\rm Pic}\le\log(\partial_F)$, so that $D'$ is close to~$D$.
 
In a similar way one can check that a
given divisor $D=(I,u)$ is reduced.  First of all we must have that $u=N(I)^{-1/n}$.
Then we check that  $1$ is contained in $I$.  To see  whether or not $1$ is a {\it minimal} element of~$I$, 
we need to make sure that the box 
$$B=\{(y_{\sigma})\in F_{\RR}:\hbox{ $|y_{\sigma}|<1$ for all~$\sigma$.}\}.
$$
contains no non-zero points of the lattice $I\subset F_{\RR}$. The box $B$ contains all vectors of length at most~1. On the other hand, every vector in~$B$ has length at most~$\sqrt{n}$.

If the first vector $\b_1$ of the LLL-reduced basis has length less than~$1/\sqrt{n}$, it is contained in~$B$ and the element $1\in I$ is {\it not} minimal. In this case we are done. 
Suppose therefore that we have $\abs{\b_1}\ge1/\sqrt{n}$.
It suffices now to compute all vectors $\bf x$ in the lattice that have length less than  
$\sqrt{n}$ and see whether they are in the box $B$ or not.
By  Corollary~10.2, the vectors ${\bf x}=\sum_{i=1}^nm_i\b_i$ of length at most~$\sqrt{n}$ have the property
that
$$
|m_{i}|\quad\le \quad 2^{{n-1}\over 2}\left({3\over 2}\right)^{n-i}{{\abs{\bf x}}\over{\abs{\b_1}}}\le  2^{{n-1}\over 2}
\left({3\over 2}\right)^{n-i}n.
$$
So, the number of vectors to be checked is bounded independently of the discriminant of~$F$.
This completes the description of the algorithm.
Both algorithms run in time polynomial in~$\log|\Delta_F|$, $\log\abs{u}$ and   the logarithmic height of~$N(I)$. 

\proclaim Algorithm 10.4.  ({\it Composition algorithm.}) Given two reduced Arakelov divisors $D=d(I)$ and~$D'=d(J)$, 
 compute a reduced divisor that is close to the sum $D+D'$ in~${\rm Pic}_F^0$.

\smallskip\noindent{\bf Description.}  One first
adds $D$ and $D'$ as divisors. Since $N(I^{-1}), N(J^{-1})\le\partial_F$, the result $(IJ, N(IJ)^{-1/n})$ can
be computed in time polynomial in~${\rm log}|\Delta_F|$. Then one reduces the result by means of
Algoithm~10.3. Since we have $N(IJ)^{-1}\le\partial_F^2$, the running time of this second step is also
polynomial in~${\rm log}|\Delta_F|$.

\proclaim Algorithm 10.5. ({\it Inversion algorithm.}) Given a reduced Arakelov divisor $D=d(I)$,  compute 
a reduced divisor that is close to  $-D$ in ${\rm Pic}_F^0$. 

\smallskip\noindent{\bf Description.}  One just
computes the inverse ideal $I^{-1}$ and reduces the divisor $d(I^{-1})$ by means of 
Algorithm~10.3. Since $N(I^{-1})\le\partial_F$, the running time of this algorithm is also
polynomial in~$\log|\Delta_F|$.
\medskip
Before describing the next algorithm, it is convenient to prove a lemma.

\proclaim Lemma 10.6. Let $D=(I,u)$ be an Arakelov divisor of degree~0 and let $\epsilon>0$.
Then every reduced divisor at distance at most~$\epsilon$ from~$D$ is of the form $d(I\mu^{-1})$
where $\mu$ is a minimal element of $I$ satisfying 
$$\abs{{\mu}}_D<\sqrt{n}\,e^{2\epsilon}\abs{y}_D,\qquad\hbox{ for all non-zero $y\in I$.}
$$
In particular, the inequality holds  for a  non-zero  $y\in I$ that is shortest with respect to the metric of~$P$.

\smallskip\noindent{\bf Proof.}  Let $D'$ be a  reduced divisor for which we have $\abs{D-D'}_{\rm Pic}<\epsilon$.
Then we have $D'=d(I\mu^{-1})$ for some minimal element $\mu\in I$. By Proposition~6.2 there is a unit $\varepsilon$ so that  for $D''=D-(\mu)+(O_F,|\varepsilon|)$ we have
$$
e^{-\abs{D'-D''}_{\rm Pic}}\quad\le\quad{{\abs{x}_D'}\over {\abs{x}_{D''}}}
\quad\le\quad e^{\abs{D'-D''}_{\rm Pic}},\qquad\hbox{for every $x\in I\mu^{-1}$.}
$$
We multiply  $\mu$ by $\varepsilon$.  Then $\mu$ remains a minimal element of~$I$ and the divisor 
$D'$ does not change.
But now   $D''$ is equal to~$D-(\mu)$.  Since we have $\abs{D-D'}_{\rm Pic}=\abs{D'-D''}_{\rm Pic}$,
the inequality above and   Proposition~7.1 imply that
$$
\abs{\mu}_D=\abs{1}_{D-(\mu)}\le e^{\epsilon}\abs{1}_{D'}\le e^{\epsilon}\sqrt{n}\abs{x}_{D'}
\le e^{2\epsilon} \sqrt{n}\abs{x}_{D-(\mu)}
=e^{2\epsilon}\sqrt{n}\abs{x\mu}_{D},
$$
for any non-zero $x\in I\mu^{-1}$. It follows that  we have
$\abs{\mu}_D\le \sqrt{n}e^{2\epsilon}\abs{y}_{D}$
for all non-zero~$y\in I$. 
This proves the proposition.

\medskip
I owe   the following algorithm to Hendrik Lenstra. See  [\BUJNT, \BUHAB, \THIE] for a different approach.

\proclaim Algorithm 10.7. ({\sl Scan algorithm.}) 
 Let $D=(I,u)$ be an Arakelov divisor of degree~0 and let $R>0$. Compute all reduced
Arakelov divisors in a ball in the Arakelov class group ${\rm Pic}_F^0$ of radius $R$ and center~$D$
in time polynomial in $\log|\Delta_F|$ and linear in the volume of the ball.

\smallskip\noindent{\bf Description.} Choose $\epsilon,\epsilon'\in\RR$ such that $0<\epsilon'<\epsilon<R$.
Inside the open ball of divisors in~${\rm Pic}_F^0$ having distance at most $R+\varepsilon$ from~$D$,
we compute a {\it web}   of regularly distributed points. 
The points $P$ in the web are at most $\varepsilon$  and at least $\varepsilon'$ apart, say.
By Theorem~7.4  every $P$ is the class of a divisor of the form $D'+(O_F,v)$ for 
some reduced divisor~$D'=d(J)$ and a totally positive  $v\in F_{\RR}^*$ satisfying
$\abs{v}_{\rm Pic}<\log(\partial_F)$. 
Moreover, LLL-reduced bases for the  lattice associated to each~$P$ can be computed in time polynomial in 
$\log|\Delta_F|$.

By Lemma~10.6, the reduced divisors we are looking for are among  the divisors 
of the form $d(J\mu^{-1})$ where $D'=d(J)$ is reduced, $P=D'+(O_F,v)$ is in the web
and $\mu\in J$ is a {\it minimal} element for which $\abs{\mu}_P$ is at most $e^{2\epsilon}\sqrt{n}$ 
times the length of a  shortest non-zero element~$x\in J$.
So, it suffices to compute the elements $\mu$ for all $P$ in the web.
For a given $P$, Corollary~10.2 says that the number of vectors $\mu\in J$ 
 of length at most $e^{2\epsilon}\sqrt{n}$ times the length of the shortest non-zero vector, 
 is  bounded independently of $P$ and even of the discriminant of~$F$.
 They can be computed in time polynomial in~$\log|\Delta_F|$.
 Minimality of the elements~$\mu$
 can be tested by means of Algorithm~10.3. Finally, since the divisors $P$ are at least $\epsilon'$ apart,
  the number of points in the web is
 proportional to the volume of the ball. This completes the description of the algorithm.
 
\proclaim Algorithm 10.8. ({\sl Jump algorithm.})  Given the coefficients of a divisor $D=\sum_{\pp}n_{\pp}\pp+\sum_{\sigma}x_{\sigma}\sigma$ of degree~0,  compute
a reduced Arakelov divisor whose image  in~${\rm Pic}_F^0$ has 
distance less than $\log(\partial_F)$ from~$D$.

\smallskip\noindent{\bf Description.} 
We assume that at most~$O(\log|\Delta_F|)$
coefficients of~$D$ are non-zero and that the coefficients themselves are of size~$|\Delta_F|^{O(1)}$. Directly
applying the reduction algorithm to~$D$ is not a very good idea, since the  LLL-algorithm and therefore
the reduction algorithm
run in time polynomial in $\log\abs{u}$ which is {\it exponential} in the terms of the $\log(x_{\sigma})$.
Therefore we proceed differently. 

It is easy to compute a reduced divisor that is close to  $d(I)$. Here $I$ is the ideal~$\prod_{\pp}\pp^{n_{\pp}}$. Indeed, for each prime ideal $\pp$ with $n_{\pp}\not=0$, we compute a reduced divisor $D'$ close to 
$d(\pp)=(\pp,N(\pp)^{-1/n})$ and then by successive squarings, multiplications and 
reductions we compute  a reduced divisor close to~$d(I)$. 
Here we use the binary expansion of the coefficient~$n_{\pp}$.
The coefficients at the primes~$\pp$ of the divisor $D-d(I)$ are all equal to zero.
Therefore it suffices to explain how to compute a reduced divisor that is close to a given 
divisor $D$  of the form
 $D=(O_F,v)$ for some~$v$.

We first compute a list of 
$r_1+r_2-1$ reduced divisors $E_i$ in  the connected component of identity of~${\rm Pic}_F^0$
that are equivalent to divisors $(O_F,v_i)$ 
that have the property that the vectors $v_i$ form
a  reasonably orthogonal basis of short vectors of the, multiplicatively written,
metric vector space~$(\prod_{\sigma}\RR^*_+)^0$. This can be accomplished  as follows. Pick
$w_1\in(\prod_{\sigma}\RR^*_+)^0$ of length ${\rm log}|\Delta_F|$ at random and reduce the divisor
$(O_F,w_1)$ to $E_1=d(f_1^{-1}O_F)$ for some $f_1\in O_F$. Then $E_1$ is equivalent to $(O_F,v_1)$ where
$v_1=|f_1|/|N(f_1)|^{1/n}w_1$. When $E_1,\ldots,E_i$ are constructed, then pick
$w_{i+1}$ of length ${\rm log}|\Delta_F|$ in the orthogonal complement of the span of
$\{v_1,\ldots,v_i\}$ inside $(\prod_{\sigma}\RR)^0$ and reduce~$(O_F,w_{i+1})$ to $E_{i+1}$, which is
equivalent to $(O_F,v_{i+1})$ where $v_{i+1}=|f_{i+1}|/|N(f_{i+1})|^{1/n}w_{i+1}$ for some $f_{i+1}\in O_F$.
This completes the description of the calculation of the `good basis'~$\{E_i\}$.

To express~$D=(O_F,v)$ as a sum of the divisors $E_i$, we 
solve a linear system and write $v=\prod_iv_i^{\lambda_i}$ and 
let $m_i$ denote the integers nearest to $\lambda_i$ for
$1\le i\le r_1+r_2-1$. Then the coefficients $m_i$ are at most~$O(|\Delta_F|^{O(1)})$.
Using the binary expansion of the coefficients $m_i$ we compute the sum $\sum_im_iE_i$
by means of successive compositions and reductions. The result is a reduced divisor that is 
quite close to $(O_F,v)$. One can get as close as
${\rm log}(\partial_F)$ to~$D$ by additional compositions with the divisors~$E_i$ or
by adjusting the infinite components and reducing. We leave this to the reader.
The amount of calculations to do all this is bounded by~${\rm log}|\Delta_F|^{O(n)}$.

This completes the description of the algorithm.
\medskip

We leave the reader the task of modifying these  algorithms so that they work
for the group~$\widetilde{\rm Div}_F^0$  of {\it oriented} divisors and for the 
oriented Arakelov class group~$\widetilde{\rm Pic}_F^0$. 
The only difference is that the unit $u$ of an extended Arakelov divisor
$D=(I,u)$ is a complex rather than a positive real number. The image of the set of reduced Arakelov
divisors in this group is probably also reasonably dense in~$\widetilde{\rm Pic}_F^0$
and that's all we need for the Jump  Algorithm to work. See Question~9.1.

\medskip\noindent
{\bf Application 10.9.} As an application we present an algorithm to compute the function
$h^0(D)$ that was introduced in~[\EF]. For an Arakelov divisor $D=(I,u)$, the number
$h^0(D)$ should be viewed as the arithmetic analogue of the dimension of the space of global sections of a divisor~$D$ on an algebraic curve. The number $h^0(D)$  depends only on the class of
$D$ in~${\rm Pic}_F^0$ and is defined as
$$
h^0(D)={\rm log}(\sum_{f\in I}\exp({-\pi\abs{f}_D^2})).
$$
See section~4 for the close relation between the function $h^0(D)$ and the Hermite
constant~$\gamma(D)$ of the ideal lattice associated to~$D$.
Since the short vectors $f\in I$ contribute the most to this exponentially quickly converging sum, 
the function $h^0(D)$ can be evaluated most efficiently when we know a good, i.e., a reasonably orthogonal
basis for $I$. As we explained above, a direct application of a lattice reduction algorithm to $D$ may
be very time consuming. Therefore we apply the {\it Jump algorithm}. We jump to a reduced divisor
$D'=d(J)$ close to $D$ in~${\rm Pic}_F^0$. Then $D$ is equivalent to $D'+(O_F,v)$ for some short
$v\in F_{\RR}^*$ and
$$\eqalign{h^0(D)&=h^0(D'+(O_F,v))={\rm log}(\sum_{f\in J}\exp({-\pi\abs{f}_{D'+(O_F,v)}^2})),\cr
&=
{\rm log}(\sum_{f\in J}\exp\left({-{{\pi}{N(J)^{-2/n}}}\sum_{\sigma}{\rm
deg}(\sigma)|\sigma(f)|^2v_{\sigma}^{2}}\right)).\cr}
$$
Since $D'$ is reduced and the vector $v=(v_{\sigma})_{\sigma}$ is short, an LLL reduced basis for the lattice associated to
$D'+(O_F,v)$ can be computed efficiently . Such a basis can be computed
efficiently since $J^{-1}$ is an integral ideal of norm at most~$|\Delta_F|^{1/2}$. This completes the
description of the algorithm to compute~$h^0(D)$.

\beginsection 11. A deterministic algorithm.

In this section we describe a {\it deterministic} algorithm to compute the Arakelov class group of a number field~$F$ of degree~$n$ and discriminant $\Delta_F$.  It runs  in time proportional to $\sqrt{|\Delta_F|}$ times a power of $\log{|\Delta_F|}$. 

\proclaim Lemma 11.1. Let $B>0$.
Then any ideal $J\subset O_F$ with
 $N(J)<B$ is of the form $J=xI^{-1}$, where
 \smallskip
 \item{--}  the Arakelov  divisor $D=(I,N(I)^{-1/n})$ is reduced;
 \item{--} the element $u=N(x)^{1/n}/|x|$  of $F_{\RR}^*$  satisfies~$\abs{u}_{\rm Pic}<\log(\partial_F)$;
 \item{--} the element  $x$  is contained in~$I$ and satisfies $\abs{x}_{D+(O_F,u)}<\sqrt{n}B^{1/n}$.
 
 \smallskip\noindent{\bf Proof.}  Suppose that $J\subset O_F$ satisfies $N(J)<B$.
 By Minkowski's Theorem there exists $y\in J^{-1}$, a shortest vector in $J^{-1}\subset F_{\RR}$ satisfying
 $|\sigma(y)|<N(J)^{-1/n}\partial_F^{1/n}$ for every~$\sigma$. We pick such an element $y$, put $x=1/y$ and $I=xJ^{-1}$.  Then the Arakelov divisor $D=(I,N(I)^{-1/n})$ is reduced. Moreover, since $xI^{-1}=J\subset O_F$, we have~$x\in I$.
 
 Writing  $u=N(x)^{1/n}/|x|$, all coordinates of the vector $N(I)^{-1/n}ux$ have absolute value
 $N(I)^{-1/n}N(x)^{1/n}=N(J)^{1/n}$ so that
 $\abs{x}_{D+(O_F,u)}=\sqrt{n}N(J)^{1/n}<\sqrt{n}B^{1/n}$.  Finally, we estimate $\abs{u}_{\rm Pic}$.
 Since $N(I)\le 1$, we have
 $$
 |u_{\sigma}|={{|N(x)|^{1/n}}\over{|\sigma(x)|}}=|\sigma(y)||N(x)|^{1/n}\le{{N(J)^{-1/n}\partial_F^{1/n}}{|N(x)|^{1/n}}}
 =N(I)^{1/n}\partial_F^{1/n}\le\partial_F^{1/n}.
 $$
 Lemma~7.5 implies than that $\abs{u}_{\rm Pic}\le\left(1-{1\over n}\right)^{1/2}\log(\partial_F)\le\log(\partial_F)$
 as required.
 
 \medskip
 It is not difficult to see that the converse of Lemma~11.1 also holds:
 any ideal $J\subset O_F$ for which the three conditions are satisfied, automatically has norm at most~$B$.

\proclaim Algorithm 11.2.  Suppose we have computed all reduced
divisors in a given connected component of the Arakelov class group~${\rm Pic}_F^0$.
In the component, detect all divisors that are  of the form $(J^{-1},N(J)^{1/n})$ with 
 $J\subset O_F$ and $N(J)<\partial_F$. 
 
\smallskip\noindent{\bf Description.}  Let $\epsilon,\epsilon'\in\RR$ such that $0<\epsilon'<\epsilon$.
 For each reduced divisor~$D=(I,N(I)^{-1/n})$ in the given connected component,  we make a web
in the ball of center~$D=(I,N(I)^{-1/n})$ and radius $\log(\partial_F)$,
 whose members $P=D+(O_F,v)=(I,N(I)^{-1/n}v)$ are at most $\epsilon$ and at least $\epsilon'$ apart. 
For each  divisor $P=D+(O_F,v)$  in  the web, 
 we compute the vectors $x$ for which we have 
 $\abs{x}_{P}\le\sqrt{n}e^{2\epsilon}\partial_F^{1/n}$.  This is done as follows.
 First we compute an LLL-reduced basis $\b_1,\ldots,\b_n$ for the lattice associated to the
 Arakelov divisor~$P$. Let $\b_1^*,\ldots,\b_n^*$ denote its Gram-Schmidt orthogonalization.
 By Lemma~10.1 we have  for any vector $x=\sum_{i=1}^nm_i\b_i$ in the lattice for which $\abs{x}_P$
 is  at most  $\sqrt{n}e^{2\epsilon}\partial_F^{1/n}$,  that
 $$
|m_i|\abs{\b_i^*}\le \left({3\over{\sqrt{2}}}\right)^{n-1}\sqrt{n}e^{2\epsilon}\partial_F^{1/n}.
 $$ 
 We simply try {\it all} coefficients  $m_i$ satisfying this inequality.

 For each such element~$x$ we then compute the corresponding ideals~$J=I^{-1}x$.
 The ideals $J$ that we compute in this way are contained in $O_F$. Moreover, every ideal 
 $J\subset O_F$ of norm at most~$\partial_F$ and for which the Arakelov divisor $(J^{-1},N(J)^{1/n})$
 lies on the given component, is obtained in this way. Indeed, if we have $N(J)<\partial_F$, Lemma~11.2 with $B=\partial_F$ implies
 that $J=xI^{-1}$ for some reduced divisor $d(I)=(I,N(I)^{-1/n})$ and some $x\in I$. Moreover, we have
 $\abs{x}_{D+(O_F,u)}<\sqrt{n}\partial_F^{1/n}$ for some  $u$ satisfying  $\abs{u}_{\rm Pic}<\log(\partial_F)$.
This means that the divisor $D+(O_F,u)$ is contained in the ball of 
center~$D=(I,N(I)^{-1/n})$ and radius $\log(\partial_F)$. Therefore there is a member $P=D+(O_F,v)$ 
of the web at distance at most $\epsilon$ from $D+(O_F,u)$. Proposition~6.2 implies then that
$$
 \abs{x}_{P}\le e^{2\epsilon}\sqrt{n}\partial_F^{1/n},
 $$
 as required.
 
 This shows that we encounter all ideals $J$ that we are after. But we'll find many more
 and we'll find each ideal many times.
  Indeed, the  vectors $x=\sum_{i=1}^nm_i\b_i$ that we consider in the computation above satisfy
  $|m_i|\abs{\b_i^*}\le \left({3\over{\sqrt{2}}}\right)^{n-1}\sqrt{n}e^{2\epsilon}\partial_F^{1/n}$ for each~$i$
  and hence
 $$
 \abs{x}_{P}\le n\left({3\over{\sqrt{2}}}\right)^{n-1}e^{2\epsilon}\partial_F^{1/n}.
 $$ It follows from the arithmetic geometric mean inequality that for the ideal $J=xI^{-1}$ we have
 $$
 N(J)=N(xI^{-1})\le n^{n/2}e^{2\epsilon n}\left({3\over{\sqrt{2}}}\right)^{n(n-1)}\partial_F.
 $$
 In order to estimate the running time of this algorithm, we estimate the number of ideals $J$ that we compute {\it and} in addition, we estimate for how many divisors $P$ in the web and how many vectors $x$, we obtain each ideal~$J$.
 By [\LEB,~Thm.6.5], the number of ideals~$J$ is bounded by $\sqrt{|\Delta_F|}$ times a power of
 $\log{|\Delta_F|}$ times a constant that depends only on the degree~$n$. Next we bound the number of times we find each ideal~$J$.
  
 First, suppose that for some divisor $P=(I,N(I)^{-1/n})+(O_F,v)$ in the web, there are two elements $x,x'\in I^{-1}$
 satisfying $\abs{x}_{P},\abs{x'}_{P}\le \sqrt{n}e^{2\epsilon}\partial_F^{1/n}$, for which the ideals $xI^{-1}$ and $x'I^{-1}$ are the {\it same}. Then we have $|\sigma(x)N(I)^{-1/n}v_{\sigma}|\le \sqrt{n}e^{2\epsilon}\partial_F^{1/n}$ for each~$\sigma$.
 Since we have $|N(v)|=1$, the product over~$\sigma$ satisfies
 $$
\prod_{\sigma}|\sigma(x)N(I)^{-1/n}v_{\sigma}|^{{\rm deg}(\sigma)}=N(xI^{-1})\ge 1.
 $$
Therefore we have
$$
-(n-1)\log(\sqrt{n}e^{2\epsilon}\partial_F^{1/n})\quad\le\quad\log|\sigma(x)N(I)^{-1/n}v_{\sigma}|
\quad\le\quad\log(\sqrt{n}e^{2\epsilon}\partial_F^{1/n})
$$
for every~$\sigma$.
We have the same inequalities for~$x'$. Therefore the unit $\varepsilon=x'/x$ satisfies
$$
-\log(\partial_F)-n\log(\sqrt{n}e^{2\epsilon})\quad\le\quad\log|\sigma(\varepsilon)|=\log\left|{{\sigma(x')}\over{\sigma(x)}}\right|
\quad\le\quad \log(\partial_F)+n\log(\sqrt{n}e^{2\epsilon}),
$$
for every~$\sigma$ and hence we have
$$
\abs{\log|\varepsilon|}\le \sqrt{n}\log(\partial_F)+n^{3/2}\log(\sqrt{n}e^{2\epsilon}).
$$
By Dobrowolski~[\DOB], there exists an absolute constant $c>0$, so that any unit $\varepsilon\in O_F^*$ that is not a root of unity satisfies $\abs{\log|\varepsilon|}>cn^{-3/2}$. Since the number of roots of unity in~$F$ is $O(n\log(n))$, the number of units satisfying the bounds above is bounded by a polynomial expression in~$\log(\partial_F)$. It follows that the number of distinct elements $x\in I$ for which the ideals $xI^{-1}$ are equal to the same ideal $J\subset O_F$ is also bounded by a polynomial expression in $\log(\partial_F)$.

  Next, suppose that an ideal  $J\subset O_F$  of norm at most~$\partial_F$ is of the form $xI^{-1}$ where
  $D=(I,N(I)^{-1/n})$  is a reduced divisor and $x\in I$ satisfies  $\abs{x}_{P}\le e^{2\epsilon}\sqrt{n}\partial_F^{1/n}$
  for some divisor  $P=D+(O_F,v)$ in the web constructed. In particular, 
 $v$ satisfies $\abs{v}_{\rm Pic}<\log(\partial_F)$.  This implies that 
 $$
\left| {{\sigma(x)}\over{N(x)^{1/n}}}v_{\sigma}^{-1}\right|=\left|\sigma(x)N(I)^{1/n}v_{\sigma}^{-1}\right|
{1\over{N(J)^{1/n}}}<{{\sqrt{n}e^{2\epsilon}\partial_F^{1/n}}\over{N(J)^{1/n}}}\le\sqrt{n}e^{2\epsilon}\partial_F^{1/n}.
$$
It follows that  the  Arakelov divisors $P$ and $(J^{-1},N(J)^{1/n})$ are rather close to one another
in~${\rm Pic}_F^0$. Indeed, we have
$$
\abs{P-(J^{-1},N(J)^{1/n})}_{\rm Pic}=\abs{(O_F,|x|N(x)^{-1/n}v_{\sigma}^{-1})}_{\rm Pic}.
$$
Since we have $\log|\sigma(x)N(x)^{-1/n}v_{\sigma}^{-1}|<\log(\sqrt{n}e^{2\epsilon}\partial_F^{1/n})$
for every infinite prime~$\sigma$, it follows from Lemma~7.5 that we have
$$
\abs{P-(J,N(J)^{1/n})}_{\rm Pic}<\log(n^{2/n}e^{2\epsilon/n}\partial_F).
$$
By Corollary~7.9, the number of reduced divisors in a ball is 
bounded by some constant, depending only on the degree
of the number field, times its volume. Therefore the number of web members $P$ for which we encounter 
a given ideal $J\subset O_F$, is bounded by a polynomial expression in~$\log(\partial_F)$.

 This completes the description and our analysis of the algorithm.

\medskip\noindent {\bf A Deterministic Algorithm.}
Finally we explain the deterministic algorithm to compute the Arakelov class group of a number field~$F$.
This algorithm seems to have been known to the experts. It was explained to me by Hendrik Lenstra.
We start at the neutral element $(O_F,1)$ of the Arakelov class group. We use Algorithm~10.3 to determine all reduced Arakelov divisors in the ball of radius $2\log(\partial_F)$ and center~$(O_F,1)$. Then we do the same with the reduced divisors $D$ we found: determine all reduced Arakelov divisors in the ball of radius $2\log(\partial_F)$ and center~$D$. Proceeding in a systematic way that is somewhat complicated to write down, we find in this way {\it all} reduced divisors in the connected component of identity.  Keeping track of their positions in terms of the coordinates in~$\prod_{\sigma}F_{\sigma}$
one computes in this way the absolute values of a set of generators of the unit group~$O_F^*$.
The running  time is proportional to the volume of the connected component of identity and is 
polynomial in~$\log{|\Delta_F|}$. 

Next we use Algorithm~11.2 and make a list $\goth L$ of  
all integral ideals $J\subset O_F$ of norm at most $\partial_F$, 
for which $(J^{-1},N(J)^{1/n})$ is on the connected component of identity. 
The amount of work is again proportional to the volume of the connected component of identity and polynomial in~$\log{|\Delta_F|}$.
By Minkowski's Theorem, the prime ideals of norm at most $\partial_F$ generate
the ideal class group of~$F$. Therefore
we check whether all prime ideals of norm at most $\partial_F$ are in the list.
This involves computing gcd's of the polynomial that defines the number field~$F$ with the polynomials $X^{p^i}-X$ for $i=1,2,\ldots,n$ for prime numbers~$p$ that are smaller than the Minkowski bound~$\partial_F$. One reads off the degrees of the prime ideals over~$p$ and hence the number of primes of norm~$p^i$ for $i=1,2,\ldots$.  The amount of work 
is linear in the length of the list and polynomial time in~$\log p$ for each prime~$p$.  If  all prime ideals of norm at most~$\partial_F$ are in the list~$\goth L$, then we are done. The class number is~1 and the Arakelov class group is connected. 

However, if we do encounter a prime number $p$, for which a prime ideal $\pp$ of norm~$p^i<\partial_F$ is missing,
then we compute it. This involves factoring a polynomial of degree~$n$ modulo~$p$. When we do this with a simple minded trial division algorithm, the amount of work is at most $p^i<\partial_F$ times a power of~$\log{|\Delta_F|}$.
By successive multiplications and reductions, we compute for $j=1,2,\ldots$
reduced divisor $D_j$ in the connected components of the Arakelov class groups  that contain 
divisors of the form~$(\pp^j,u)$ for some~$u$. Each time we check whether $D_j$ is already in the list~$\goth L$.
If it is, we stop computing divisors~$D_j$.

Then we repeat the algorithm, but this time we work with the connected components of the divisors ~$D_j$ rather than~$(O_F,1)$: we use Algorithm~10.3 to determine all reduced Arakelov divisors in the balls of radius $2\log(\partial_F)$ and center~$D_j$. Then we do the same with the reduced divisors we found, and so on $\ldots$. Once we have computed all reduced divisors on the connected components of~$D_j$, 
we use Algorithm~11.2 to compute all integral ideals $J\subset O_F$ of norm at most $\partial_F$, 
for which $(J^{-1},N(J)^{1/n})$ is on the connected components of the divisors~$D_j$ and we add these to the list~$\goth L$.

When we are done with this,  the list~$\goth L$ contains all integral ideals $J\subset O_F$ of norm at most $\partial_F$,
whose classes are in the group generated by the ideal class of~$\pp$. We check again whether 
all prime ideals of norm at most $\partial_F$ are in the list. If this turns out to be the case, we are done.
The ideal class group is cyclic, generated by the class of~$\pp$.  If, on the other hand, we do encounter a second prime number $q$, for which a prime ideal $\goth q$ of norm~$q^i<\partial_F$ is missing, then we compute it. We compute reduced divisors that are in the components of the powers of~$\goth q$ $\ldots$ etc.

For each new prime that we find is {\it not} in the list~$\goth L$, we factor a polynomial and the amount of work to do this is at most~$\partial_F$. However, since the ideal class group has order at most $\sqrt{|\Delta_F|}$ times power of 
$\log|\Delta_F|$, we need to do this at most $\log|\Delta_F|$ times. As a result this algorithm takes time
at most $\sqrt{|\Delta_F|}$ times power of ~$\log|\Delta_F|$.

\beginsection 12. Buchmann's algorithm.

In this section we briefly sketch Buchmann's algorithm~[\BUDPP,~\BUDU] for computing the Arakelov divisor
class group and, as a corollary, the class group and regulator of a number field~$F$. This algorithm
combines the infrastructure idea with an algorithm for complex quadratic number fields presented by
J.~Hafner and K.~McCurley~[\MCCUR] in~1989. When we fix the degree of~$F$, the algorithm is
under reasonable assumptions  subexponential in the discriminant of the number field~$F$. A practical
approach is described in~[\CO, section~6.5]. The algorithm has been implemented 
in the LiDIA, MAGMA and PARI 
software packages~[\LID, \MAG, \PARI]. See also~[\THIE].

Let $F$ be a number field of degree~$n$. The structure of Buchmann's algorithm is very simple. 
Our first description involves the Arakelov class group ${\rm Pic}_F^0$ rather than the
oriented group~$\widetilde{\rm Pic}_F^0$.

\smallskip\noindent {\bf Step 1.} {\sl Estimate the volume of~${\rm Pic}_F^0$.}
By Prop.~6.5 the volume of the compact Lie group ${\rm Pic}_F^0$ is given by
$$
{\rm vol}({\rm Pic}^0_F)=
{{w_F\,\sqrt{n}}\over{2^{r_1}(2\pi\sqrt{2})^{r_2}}}\cdot|\Delta_F|^{1/2}\cdot
\mathop{\rm
Res}\limits_{s=1}\zeta_F^{}(s).
$$
The computation of $r_1$, $r_2$ and $w_F=\#\mu_F$ is easy. The discriminant is computed as a byproduct
of the calculation of the ring of integers~$O_F$. Approximating the residue of 
the zeta function 
$$\zeta_F^{}(s)\,\,=\,\,\prod_{\goth p}\left(1-{1\over{N({\goth p})^{s}}}\right)^{-1}
$$
in $s=1$ is done by dividing $\zeta_F(s)$ by the  zeta function of~$\QQ$
and by directly evaluating a truncated Euler product
$$
\prod_{p\le X}\left(1-{1\over p}\right)\prod_{\pp|p}\left(1-{1\over{N(\pp)}}\right)^{-1}.
$$
This involves factoring the ideals $pO_F$ for all prime numbers $p<X$. See~[\CO] for efficient
methods to do this. The Euler product converges rather slowly.
Under assumption of the Generalized Riemann Hypothesis for the
zeta function of~$F$,  using the primes $p<X$, the relative error is
$O(X^{-1/2}{\rm log}|\Delta_FX|)$.  Here the O-symbol only depends on
the degree of the number field~$F$. See~[\BUW,~\SMC].
Therefore, there is a constant $c$ only depending on the degree of~$F$, 
so that if we truncate the Euler product at $X= c\,{\rm log}^2|\Delta_F|$, 
the relative error in the approximation of ${\rm vol}({\rm Pic}_F^0)$
is at most~$1/2$.

\smallskip\noindent {\bf Step 2.} {\sl Compute a factor basis.}  We compute a factor base ${\cal
B}$ i.e., a list of prime ideals $\pp$ of~$O_F$ of norm less than $Y$ for some~$Y>0$. Computing a
factor basis involves factoring the ideals $pO_F$ for various prime numbers~$p$. It  is
convenient to do this alongside the computation of the Euler factors in Step~1.  We add the
infinite primes to our factor basis. 
By normalizing, we obtain in this way a factor basis of Arakelov divisors of degree~0.
The factor basis should be so large that
the natural homomorphism
$$
\left(\mathop{\oplus}\limits_{\pp\in{\cal B}}\ZZ\times\mathop{\oplus}\limits_{\sigma}\RR\right)^0\quad\longrightarrow\quad{\rm Pic}_F^0
$$
is surjective. By Prop.~2.2 this means that the classes of the primes in~${\cal B}$
must generate the ideal class group. Under assumption of the
Generalized Riemann  Hypothesis for the $L$-functions $L(s,\chi)$ associated to characters $\chi$
of the ideal class group $Cl_F$ of~$F$, this is the case for
$Y>c'\log^2|\Delta_F|$ for some constant $c'>0$ that only depends on the degree of~$F$. 
Taking $\cal B$ this big, we have that
$$
{\rm Pic}_F^0\quad=\quad\left(\mathop{\oplus}\limits_{\pp\in{\cal B}}\ZZ\times\mathop{\oplus}\limits_{\sigma}\RR\right)^0/H
$$
where $H$ is the discrete subgroup of principal divisors of
${\cal B}$-units, i.e., the group $H$ consists of divisors $(f)$ where $f\in F^*$ 
are elements whose prime factorizations involve only prime ideals~$\pp\in{\cal B}$.

\smallskip\noindent {\bf Step 3.} {\sl Compute many elements in~$H$.}
 An Arakelov divisor $D=(I,u)$ is called $\cal B$-smooth if $I$ 
is a product of powers of primes in~$\cal B$.  We need to find elements $f\in F^*$ for which~$(f)$ is $\cal B$-smooth
and hence~$(f)\in H$. This is achieved by
repeatedly doing the following.
For at most $O({\rm log}|\Delta_F|)$ prime ideals $\pp\in{\cal B}$ pick random exponents $m_{\pp}\in\ZZ$ 
of absolute value not larger than $|\Delta_F|$. In addition, pick random
$x_{\sigma}\in\RR$ of absolute value not larger than~$|\Delta_F|$. Replacing $x_{\sigma}$ by
$x_{\sigma}N(D)^{-1/n}$, scale  the Arakelov divisor
$D=\sum_{\pp}m_{\pp}\pp+\sum_{\sigma}x_{\sigma}\sigma$ so that it acquires degree zero. Then the class of $D$ is a random  element of~${\rm Pic}^0_F$. We use the {\sl Jump Algorithm} described in section~10 and `jump to~$D$'. The result is a reduced divisor $D'=(I,N(I)^{-1/n})$ whose image in ${\rm Pic}_F^0$ is not too far from the image of~$D$.
This means that
$$
D=(f)+D'+(O_F,v)
$$
for some $f\in F^*$ and  $v=(v_{\sigma})\in \left(\prod_{\sigma}\RR^*_+\right)^0$ for which $\abs{v}_{\rm Pic}$ is small,
at most $\log(\partial_F)$ say.
There is no need to compute $f$, but when one applies the Jump Algorithm one should keep
track of the infinite components and compute~$v$ or its logarithm.

Since the divisor $D$ is random, it seems reasonable to think of
 the reduced divisor $D'=(I,N(I)^{-1/n})$ as being  `random'  as well.
Next we attempt to factor the integral ideal $I^{-1}$ into a product of prime ideals $\pp\in{\cal B}$.
Since $D'$ is random and since the norm of $I^{-1}$  is at most 
$\partial_F=\hbox{$\left({2\over{\pi}}\right)^{r_2}$}|\Delta_F|^{1/2}$ and hence  relatively small, we have
a fair chance to succeed. If we do, then we have 
$D'=\sum_{\pp\in{\cal B}}n_{\pp}\pp+\sum_{\sigma}y_{\sigma}\sigma$ and hence~$(f)\in H$.
This factorization leads to a relation of the form
$$
(f)=D-D'-(O_F,v)=
\sum_{\pp\in{\cal B}}(m_{\pp}-n_{\pp})\pp+\sum_{\sigma}(x_{\sigma}-y_{\sigma}+v_{\sigma})\sigma.
$$
In this way we have computed an explicit element in~$H$.

Since we want to find many such relations, we need to be successful relatively often.
In other words, the `random' reduced divisors $D'$ that we obtain, should be {\it $\cal B$-smooth} relatively often.
This is the weakest point of our analysis of the algorithm. In section~9 the set
${\rm Red}_F''$  of Arakelov divisors $d(I)$ for which $1\in I$ is
primitive and $N(I^{-1})\le\sqrt{|\Delta_F|}$ was introduced.
Under the assumption of the Generalized Riemann  Hypothesis, Buchmann and Hollinger~[\BUHO] showed
that when $Y\approx\exp(\sqrt{\log|\Delta_F|})$, the proportion of ${\cal B}$-smooth ideals $J$
with $d(J^{-1})\in{\rm Red}''_F$ is at least
$\exp(-\sqrt{\log|\Delta_F|}\,\log\log|\Delta_F|)$. Here the Riemann  Hypothesis for the zeta-function
of the normal closure of~$F$ is used to guarantee the existence of sufficiently many prime ideals of norm at
most $\sqrt{|\Delta_F|}$ and degree~1. It is likely, but at present not known whether the proportion of 
${\cal B}$-smooth ideals $I$ for which $d(I)$ is contained in the subset ${\rm Red}_F$ rather than ${\rm Red}'_F$,  is {\it also} at
least $\exp(-\sqrt{\log|\Delta_F|}\,\log\log|\Delta_F|)$. 
Even if this were the case, there is the problem that the divisor
$D'$ that comes out of the reduction algorithm is  not a `random' reduced divisor. Indeed,
Example~9.5 provides examples of reduced divisors that are not the reduction of {\it any}
Arakelov divisor. These reduced divisors will never show up in our calculations, since everything
we compute is a result of the reduction algorithm.
It would be of interest to know how many such reduced divisors there may be.

For the next step we need to have computed approximately as many elements in~$H$ as the size of the factor base~$\cal B$. This implies that we expect to have to repeat the computation explained above
about $\exp(\sqrt{\log|\Delta_F|}\,\log\log|\Delta_F|)$ times.
When the discriminant $|\Delta_F|$ is large, this is more work than we need to do in Steps~1, 2 and~4.
Step 3 is in practice the dominating part of the
algorithm. It follows that the algorithm is subexponential and runs in time~$O(\exp(\sqrt{\log|\Delta_F|}\,\log\log|\Delta_F|))$.

\smallskip\noindent {\bf Step 4.} {\sl Verify that the elements computed in Step~3 actually generate~$H$.}
Let $H'$ denote the subgroup of~$H$ generated by the divisors
$(f)=\sum_{\pp\in {\cal B}}k_{\pp}\pp+\sum_{\sigma}y_{\sigma}\sigma$ that we computed in Step~3.  The quotient group
$\left(\oplus_{\pp\in{\cal B}}\ZZ\times\oplus_{\sigma}\RR\right)^0/H'$ admits a natural
map onto ~${\rm Pic}_F^0$. Its  volume is equal to the determinant of a square  matrix of
size $\#{\cal B}$ whose rows are  the coefficients of a set of $\#{\cal B}$ independent principal divisors that generate~$H'$. If the quotient of the  volume by  the estimate of ${\rm vol}({\rm Pic}_F^0)$ computed in Step~1, is less than $1/2$, then we have $H'=H$ and the group $\left(\oplus_{\pp\in{\cal B}}\ZZ\times\oplus_{\sigma}\RR\right)^0/H'$
is actually {\it isomorphic} to~${\rm Pic}_F^0$ and we are done.

In practice this means that once we have computed  somewhat more divisors
$(f)$ in~$H$ than $\#{\cal B}$, we ``reduce" the coefficient matrix.  
From the ``reduced" matrix we can read off the structure of the ideal class group
as well approximations to the logarithms of the absolute values of  a set of units $\varepsilon$ that generate
the unit group~$O_F^*$. This enables us to compute the regulator~$R_F$.

This completes our description of Buchmann's algorithm.
\medskip
It seems difficult to compute approximations to the numbers $\sigma(\varepsilon)$
themselves from approximations to their absolute values $|\sigma(\varepsilon)|$. If one wants to 
obtain such approximations, one should apply the algorithm 
above to the {\it oriented} Arakelov class group. The computations are the same, but 
rather than real, one carries complex coordinates $x_{\sigma}$ along. 
More precisely, we have that
$$
\widetilde{\rm Pic}_F^0\quad=\quad\left(\mathop{\oplus}\limits_{\pp\in{\cal B}}\ZZ\times\mathop{\oplus}\limits_{\sigma}F^*_{\sigma}\right)^0/\widetilde{H}
$$
for the discrete subgroup~$\widetilde{H}$ that consists of elements $f\in F^*$ whose prime factorizations involve only prime ideals~$\pp\in{\cal B}$. In this way one obtains approximations to $\sigma(\varepsilon_i)$ for 
a basis~$\varepsilon_i$ of the unit group~$O_F^*$.
In principle, once one has such approximations one may solve the linear system
$\sigma(\varepsilon_i)=\sum_j\lambda_{ij}\sigma(\omega_j)$  and compute $\lambda_{ij}\in\ZZ$
so that $\varepsilon_i=\sum_j\lambda_{ij}\omega_j$ for $1\le i\le r_1+r_2-1$.
However, it is well known that the size of the coefficients $\lambda_{ij}$ may grow doubly exponentially quickly in~${\rm
log}|\Delta_F|$ and  it is therefore not reasonable to ask for an efficient algorithm that computes a set of
generators of the unit group as linear combination of the basis $\omega_k$ of the additive group~$O_F$. 

What can be done efficiently, is to compute a so-called {\it compact representation} 
of a set of generators of the unit group~$O_F^*$. Briefly, this works as follows. Using the notation used in the description of
the Jump Algorithm of section~10, one finds for each fundamental unit~$\varepsilon_j$  integers $m_{ij}$ so that
$\prod_iv_i^{m_{ij}}$ is close to~$\varepsilon_j$.  The Arakelov divisors $(O_F,v_i)$ are equivalent to reduced
divisors $d(f_i^{-1})$. While jumping towards the fundamental unit,  one keeps track of the principal ideals that are
encountered on the way. For instance, if in the process one computes the sum of the divisors $(O_F,v_i)$ and
$(O_F,v_j)$ and reduces the result by means of a shortest vector $f$, then the result is equivalent to the reduced
divisor $d((ff_if_j)^{-1})$. The size of the elements $f_i$, $f_j$ and~$f$ $\ldots$ etc. is bounded by $({\rm
log}|\Delta_F|)^{O(1)}$. With a good strategy one can jump reasonably close to the unit. The number of jumps we need 
to reach this point is also bounded by~$({\rm log}|\Delta_F|)^{O(1)}$. Using the approximations to the fundamental
units and to the vectors $f_i$, $f_j$, $f$ $\ldots$ etc, we can approximate a small element $g\in F^*$, so that the
difference between the divisor we  jumped to and the fundamental unit is equivalent to a divisor of the form
$(O_F,g)$. Since $g$ is small, we can compute it in time bounded by ${\rm log}|\Delta_F|^{O(1)}$ from its the
approximations of the various $\sigma(g)$. From this we easily obtain the fundamental unit~$\varepsilon_j$.

\bigskip\bigskip\bigskip

\bibliography

\item{[\EVA]} Bayer, E.: Lattices and number fields, 
{\sl Contemp. Math.} {\bf 241} (1999), 69--84.
\item{[\BUJNT]} Buchmann, J.: On the computation of units and class numbers by a generalization
of Langrange's algorithm, {\sl J. of Number Theory} {\bf 26} (1987), 8--30.
\item{[\BUJNTT]} Buchmann, J.: On the period length of the generalized Lagrange algorithm, 
{\sl J. of Number Theory} {\bf 26} (1987), 31--37.
\item{[\BUHAB]} Buchmann, J.: {\sl Zur Komplexit\"at der Berechnung von Einheiten und Klassenzahlen
algebraischer Zahlk\"orper}, Habilitationsschrift, D\"usseldorf 1987.
\item{[\BUWINF]} Buchmann, J. and Williams, H.C.: On the infrastructure of the 
principal ideal class of an algebraic number field of unit rank one, {\sl Math. Comp.} {\bf 50}
(1988), 569--579.
\item{[\BUW]} Buchmann, J. and Williams, H.C.: On the computation of the class number
of an algebraic number field, {\sl Math. Comp.} {\bf 53} (1989), 679--688.
\item{[\BUDPP]} Buchmann, J.: {\sl A subexponential algorithm for the determination
of class groups and regulators of algebraic number fields}, pp. 27--41 in C.~Goldstein (ed):
``S\'eminaire de Th\'eorie des Nombres, Paris 1988--1989", Birkh\"auser Boston 1990.
\item{[\BUDU]}  Buchmann, J. and D\"ullmann, S.: A probabilistic class group and regulator algorithm
and its implementation, p.53--72 in Peth\"o et al (eds.): ``Computational  Number Theory"
Proceedings of the Colloquium at Debrecen 1989, De Gruyter Berlin 1991.
\item{[\BUHO]} Buchmann, J. and Hollinger, C.: On smooth ideals in number fields, {\sl J. of Number
Theory} {\bf 59} (1996), 82--87.
\item{[\CO]} Cohen, H.: {\sl A Course in Computational Algebraic Number Theory},  Graduate 
Texts in Mathematics {\bf 138}, Springer-Verlag, Berlin 1993.
\item{[\CDDO]} Cohen, H., Diaz y Diaz, F. and Olivier, M.:  Subexponential algorithm for
class group and unit computations, {\sl J. Symbolic Computation} {\bf 24} (1997), 433--441.
\item{[\DOB]} Dobrowolski, E.: On a question of Lehmer and the number of irreducible factors
of a polynomial, {\sl Acta Arithmetica} {\bf 34} (1979), 391--401.
\item{[\GR]} Groenewegen, R.P.: The size function for number fields, Proceedings of the
XXI Journ\'ees Arithm\'etiques, {\sl Journal de Th\'eorie de Nombres de Bordeaux} {\bf 13}
(2001), 143--156.
\item{[\MCCUR]} Hafner, J. and  McCurley, K.: A rigorous subexponential algorithm for computation of
class groups, {\sl Journal of the AMS} {\bf 2} (1989), 837--850.
\item{[\JU]} J\"untgen, M.: {\sl Berechnung von Einheiten in algebraischen Zahlk\"orpern
mittels des verallgemeinerten Lagrangeschen Kettenbruchalgorithmus}, Diplomarbeit D\"usseldorf 1990.
\item{[\LLL]} Lenstra, A.K., Lenstra, H.W.~and Lov\'asz, L.: Factoring polynomials with
rational coefficients, {\sl Math. Annalen} {\bf 261} (1982), 515--534.
\item{[\LEX]} Lenstra, H.W.: On the computation of regulators and class
numbers of quadratic fields. In {\sl Proc.\ Journ\'ees Arithm\'etiques  Exeter
1980}, London Math.\ Soc.\ Lect.\ Notes {\bf 56} (1982), p.\ 123--150.
\item{[\LEB]} Lenstra, H.W.: Algorithms in algebraic number theory, {\sl Bulletin of the AMS} 
{\bf 26} (1992), 211--244.
\item{[\LEM]} Lenstra, H.W.: {\sl Solving the Pell equation}, Computational Number Theory
Workshop, MSRI 2000, these proceedings.
\item{[\LENSH]} Lenstra, H.W.: {\sl Lattices}, Computational Number Theory
Workshop, MSRI 2000, these proceedings.
\item{[\LID]} LiDIA, A C++ Library For Computational Number Theory,  
Homepage: {\tt www.informatik. tu-darmstadt.de/TI/LiDIA}
\item{[\MAG]} The Magma Computational Algebra System  for   Algebra, Number Theory and Geometry,
Homepage: {\tt magma.maths.usyd.edu.au/magma}
\item{[\MAR]} Marcus, D.A.: {\sl Number Fields}, Universitext, Springer-Verlag, New York 1977.
\item{[\PARI]} Pari-GP, Homepage: {\tt www.parigp-home.de}
\item{[\SMC]} Schoof, R.: {\sl Quadratic fields and factorization}, pp.~235--286 in: H.W.
Lenstra jr. and R. Tijdeman (eds.): ``Computational Methods in
Number Theory", MC-Tracts {\bf 154-155}, Amsterdam 1982.
\item{[\SH]} Shanks, D.: The infrastructure of a real quadratic field and its applications, {\sl
Proceedings of the 1972 Number Theory Conference}, Boulder (1972) 217--224.
\item{[\SHW]} Shanks, D.: A survey of quadratic, cubic and quartic algebraic number fields (from a
computational point of view), pp. 15--40 in: {\sl Congressus Numerantium} {\bf 17} (Proc. 7th S-E
Conf. Combinatorics, graph theory and computing, Baton Rouge 1976), Utilitas Mathematica, Winnipeg
1976.
\item{[\SZ]} Szpiro, L.: Pr\'esentation de la th\'eorie d'Arakelov, p.~279--293 in ``Current
Trends in Arithmetical Algebraic Geometry", {\sl Contemporary Mathematics \bf 67},
AMS, Providence RI 1985.
\item{[\SZAST]} Szpiro, L.: Degr\'es, intersections, hauteurs, p.~11--28 in ``S\'eminaire sur
les pinceaux arithm\'e\-ti\-ques: La conjecture de Mordell", {\sl Ast\'erisque \bf
127} (1985).
\item{[\THIE]} Thiel, C.: {\sl On the complexity of some problems in algorithmic algebraic number theory},
PhD Thesis, Universit\"at des Saarlandes, Saarbr\"ucken 1995.
\item{[\EF]} Van der Geer, G. and Schoof, R.:  Effectivity of Arakelov divisors and
the Theta divisor of a number field, {\sl Selecta Mathematica}, New Ser. 6
(2000), 377--398. Preprint 9802121 at: {\tt http://xxx.lanl.gov/list/math.AG/9802}.
\item{[\WS]} Williams, H.C. and Shanks, D.: A note on class number one in pure cubic fields, {\sl
Math. Comp.} {\bf 33} (1979), 1317--1320.
\item{[\WDS]} Williams, H.C., Dueck, G. and Schmid, B.: A  rapid method of evaluating the regulator
and class number of a pure cubic field, {\sl Math. Comp.} {\bf 41} (1983), 235--286.

\bye